%% file: multifractal-thermo.tex
\documentclass[11pt]{amsart}

\usepackage{amssymb}
\usepackage{enumerate}
\usepackage{graphicx}

\include{defs}

\newcommand{\PP}{\textbf{(P)}}

\newcommand{\tG}{\tilde{G}}

\newcommand{\LLL}{\mathcal{L}}
\newcommand{\KKK}{\mathcal{K}}
\newcommand{\DDD}{\mathcal{D}}
\newcommand{\PPP}{\mathcal{P}}
\newcommand{\EEE}{\mathcal{E}}
\newcommand{\FFF}{\mathcal{F}}
\newcommand{\ZZZ}{\mathcal{Z}}

\newcommand{\Mf}{\MMM^f}
\newcommand{\Mfe}{\MMM^f_E}
\newcommand{\Mfa}{\MMM^f_\alpha}

\DeclareMathOperator{\spn}{span}
\DeclareMathOperator{\inter}{int}
\DeclareMathOperator{\Capac}{Cap}

\newcommand{\llambda}{\underline{\lambda}}
\newcommand{\ulambda}{\overline{\lambda}}

\newcommand{\lhtop}{\underline{Ch}_\mathrm{top}}

\newcommand{\lCap}{\underline{\Capac}}
\newcommand{\uCap}{\overline{\Capac}}

\newcommand{\lF}{\underline{\FFF}}
\newcommand{\uF}{\overline{\FFF}}

\theoremstyle{plain}
\newtheorem{theorem}{Theorem}[section]
\newtheorem{proposition}[theorem]{Proposition}
\newtheorem{corollary}[theorem]{Corollary}
\newtheorem{lemma}[theorem]{Lemma}
\newtheorem{thma}{Theorem}
\newtheorem*{thma*}{Theorem}

\theoremstyle{definition}
\newtheorem{definition}{Definition}[section]

\theoremstyle{remark}
\newtheorem{example}[theorem]{Example}
\newtheorem*{remark}{Remark}
\newtheorem{problem}[theorem]{Problem}

\numberwithin{equation}{section}

\title{The thermodynamic approach to multifractal analysis}
\author{Vaughn Climenhaga}
\address{Department of Mathematics \\ University of Maryland \\ College Park, MD 20742, USA.}
\email{climenhaga@math.umd.edu}
\urladdr{http://www.math.umd.edu/$\sim$climenhaga/}

\begin{document}

\date{\today}
\begin{abstract}
Most results in multifractal analysis are obtained using either a thermodynamic approach based on existence and uniqueness of equilibrium states or a saturation approach based on some version of the specification property.  A general framework incorporating the most important multifractal spectra was introduced by Barreira and Saussol, who used the thermodynamic approach to establish the multifractal formalism in the uniformly hyperbolic setting, unifying many existing results.  We extend this framework to apply to a broad class of non-uniformly hyperbolic systems, including examples with phase transitions, and obtain new results for a number of examples that have already been studied using the saturation approach.  We compare the thermodynamic and saturation approaches and give a survey of many of the multifractal results in the literature.
\end{abstract}

\thanks{This work is partially supported by NSF grant 0754911.  A preliminary announcement (without proofs) of some of the results in this paper appeared in \emph{Electronic Research Announcements}~\cite{vC10}.  A number of these results (and their proofs) have been made available in the permanent preprint~\cite{vC10c}.}

\maketitle

\section{Introduction}

\subsection{Overview and basic concepts}

Many important characteristics of dynamical systems are given as local asymptotic quantities, such as Birkhoff averages, Lyapunov exponents, local entropies, and pointwise dimensions, which reveal information about a single point or trajectory.  It is of interest to study the level sets for these quantities:  when this is done using measure theory, one is led into ergodic theory; when it is done using dimension theory, one is led into multifractal analysis.

The study of multifractal analysis is multifaceted, and our purpose in this paper is two-fold.  First, we consider multifractals from the point of view set out in~\cite{BPS97} and prove results generalising those in~\cite{BS01}; in particular, we give purely thermodynamic conditions for the multifractal formalism to hold, allowing multifractal analysis to be carried out for very general systems by studying the topological pressure.  Broadly speaking, one may say that the results in~\cite{BS01} unify the menagerie of multifractal spectra in the uniformly hyperbolic case, and the results of the present paper extend this unification to many non-uniformly hyperbolic systems.

Secondly, we give an account of the relationship between this thermodynamic approach to multifractals and other approaches in the literature, particularly the saturation approach, with the aim of elucidating the similarities and distinctions between these works.

The basic scheme of multifractal analysis is as follows.  Given a dynamical system $f\colon X\to X$, we begin with some local asymptotic quantity $\phi_\infty(x) := \lim_{n\to\infty} \phi_n(x)$, where $\phi_n(x)$ is some sequence of functions and the limit is not required to exist everywhere.  Typically each of the functions $\phi_n$ is continuous, but we should \emph{not} expect $\phi_\infty$ to depend continuously on $x$.  Rather, we should expect the level sets $K_\alpha := \{x\in X \mid \phi_\infty(x)=\alpha\}$ to be dense in $X$ for every value of $\alpha$ where they are non-empty.  Given this \emph{multifractal decomposition}, we define a \emph{multifractal spectrum} $\FFF(\alpha) := \dim K_\alpha$, where $\dim$ is a Carath\'eodory dimension characteristic~\cite{yP98} such as Hausdorff dimension or topological entropy.

Given the manner in which $\FFF(\alpha)$ is defined, it is \emph{a priori} rather surprising to find that $\FFF$ is in many cases a real analytic and concave function of $\alpha$.  This so-called \emph{multifractal miracle} is closely related to the fact that the function $\FFF$ is often determined by the topological pressure function $P\colon C(X) \to \RR$, which has been well-studied as a key component of the thermodynamic formalism; in particular, $P$ is well-known to be convex and is often analytic on certain subspaces of $C(X)$.

One of the central goals of multifractal analysis is to determine in exactly which cases the expected relationship between topological pressure and multifractal spectra goes through: where we do not give precise results, we will refer to this state of affairs by saying ``the formalism holds''.

The multifractal spectrum $\FFF$ depends on both the local and global quantities, and so by varying these, we may consider many different spectra.  In this paper, we follow Barreira and Saussol~\cite{BS01} and consider spectra defined by functions $\ph,\psi, u \in C(X)$.  The local quantity is given by $\phi_n(x) = S_n\ph(x)/S_n\psi(x)$, where $S_n$ denotes the $n$th Birkhoff sum, $S_n\phi(x) = \sum_{k=0}^{n-1} \phi(f^k(x))$.  The global quantity $\dim_u$ is a dimensional quantity depending on the function $u$ in a manner that will be made precise later on in Section~\ref{sec:pressure}; its key property is that by an appropriate choice of $u$, it can be made equal to either the topological entropy or the Hausdorff dimension, both of which are commonly used to construct multifractal spectra.

A survey describing various well-studied multifractal spectra can be found in~\cite{BPS97}.  In many important cases, these spectra fit into the above framework via appropriate choices of $\ph$, $\psi$, and $u$.  We describe these briefly now; precise statements will be given in Sections~\ref{sec:coarse}--\ref{sec:dimspec}.

\begin{enumerate}
\item If $u\equiv 1$ and $\psi\equiv 1$, then the $u$-dimension reduces to the topological entropy, while the ratio level sets $K_\alpha$ become the level sets for Birkhoff averages of $\ph$.  Thus $\FFF(\alpha)=\BBB_E(\alpha)$ is the entropy spectrum for Birkhoff averages.
\item If $\mu$ is a Gibbs measure for $\ph$ and $u\equiv \psi \equiv 1$, then the Birkhoff averages of $\ph$ and the local entropies of $\mu$ have the same level sets, so $\FFF$ determines the entropy spectrum for local entropies $\EEE_E$.
\item If $f$ is conformal ($Df$ is a scalar multiple of an isometry), $\mu$ is a Gibbs measure for $\phi$, and $\psi=\log \|Df\|$, then writing $\ph=P(\phi)-\phi$, the ratios $S_n\ph/S_n\psi$ converge to the pointwise dimension of $\mu$, so the ratio level sets $K_\alpha$ become the level sets for pointwise dimensions of $\mu$.  If in addition $u=\log \|Df\|$, then the $u$-dimension becomes the Hausdorff dimension and $\FFF(\alpha)=\DDD_D(\alpha)$ is the dimension spectrum for pointwise dimensions.
\item If $f$ is conformal, $\ph=\log \|Df\|$, and $\psi\equiv 1$, then the ratio level sets become the level sets for Lyapunov exponents.  If $u\equiv 1$, then $\FFF(\alpha)$ is the entropy spectrum for Lyapunov exponents $\LLL_E$; if $u=\log \|Df\|$, then it is the dimension spectrum for Lyapunov exponents $\LLL_D$.
\end{enumerate}

When $u=\psi$, as in most of the above examples, the spectrum $\FFF(\alpha)$ is the Legendre transform of a certain thermodynamic function.  For $u=\ph$, the spectrum itself is not a Legendre transform, but is closely related to one.  For the more general case of \emph{mixed multifractal spectra} where $u$ differs from both $\ph$ and $\psi$, Legendre transforms cannot be used directly, but we can still obtain the spectrum from thermodynamic considerations.

All this has been thoroughly discussed in~\cite{BS01}:  let $X$ be a compact metric space, $f\colon X\to X$ a continuous map, and $\ph, \psi, u$ continuous functions with $\psi,u>0$.  Denote by $\Mf(X)$ the space of $f$-invariant Borel probability measures on $X$.  Suppose that the entropy map $h\colon \Mf(X)\to \RR$ is upper semi-continuous and the pressure function $P\colon C(X)\to\RR$ is differentiable on the span of $\{\ph,\psi,u\}$; then it is shown in~\cite{BS01} that the multifractal formalism holds for $\FFF = \FFF_u^{\ph,\psi}$.  One goal of this paper is to show that these results (or suitably weakened versions) hold in a more general setting:
\begin{enumerate}
\item the requirement that $\psi$ and $u$ be uniformly positive is replaced by weaker asymptotic positivity criteria, which in particular allow the possibility that $u=0$ at a fixed point of $f$---see Section~\ref{sec:MP};
\item the pressure function need not be differentiable on the entire subspace spanned by $\{\ph,\psi,u\}$, but only on a subset thereof (with a concomitant weakening of the results)---see Theorem~\ref{thm:full-measure};
\item even if the pressure function is not differentiable (there is a phase transition), thermodynamic considerations can still sometimes be used to guarantee that the formalism holds---see Corollary~\ref{cor:phase}.
\end{enumerate}

Besides the spectra described in~\cite{BPS97} and the mixed spectra studied in~\cite{BS01}, there have been many related objects studied under the name ``multifractal analysis''.  Before moving on, we pause to mention a few of these; we will examine them in greater depth later on.

\begin{enumerate}
\item The dimension spectrum for pointwise dimensions $\DDD_D$ is defined in terms of a measure $\mu$ on a separable metric space $X$ and does not require a dynamical system for its definition.  The case where $\mu$ has some self-similarity arising from a geometric construction has been widely studied (see Section~\ref{sec:DD}).
\item The local quantity $\phi_\infty(x)$ is defined as a limit, but there are typically many points $x$ at which this limit does not exist.  The set of such points is disjoint from all level sets and is generally invisible to invariant measures, but may have large dimension.  In recent years, the study of this \emph{irregular set} has grown significantly (see Section~\ref{sec:irreg}).
\item The spectra defined in terms of level sets are sometimes referred to as \emph{fine multifractal spectra}.  For each such spectrum there is an associated \emph{coarse multifractal spectrum} (see Sections~\ref{sec:coarse} and~\ref{sec:coarse2}).
\item One may consider local quantities $\phi_n(x)$ that are not given as ratios of Birkhoff sums, or that take values in a more general space than $\RR$.  For example, by considering $\phi_n$ with values in $\RR^d$, one can study simultaneous level sets for Birkhoff averages of $d$ functions (see Section~\ref{sec:multidim}).
\end{enumerate}

\subsection{Description of approach}

The multifractal spectrum that has been the most extensively studied is the dimension spectrum for pointwise dimensions $\DDD_D$.  This spectrum was introduced in~\cite{HJKPS86}, where it is obtained for certain systems as the Legendre transform of a function $T(q)$ related to the Hentschel--Procaccia spectrum:
\begin{equation}\label{eqn:Lt1}
\DDD_D(\alpha) = T^L(\alpha) := \inf_{q\in \RR} (T(q) + q\alpha).
\end{equation}
It follows from general properties of Legendre transforms~\cite{rR70} that in this case $\DDD_D$ is concave.  Furthermore, it is true in many cases that $T$ is an analytic and strictly convex function of $q$, which implies that $\DDD_D$ is an analytic and strictly concave function of $\alpha$.

For more general spectra $\FFF(\alpha)$, a similar approach is used.  A certain family of functions $T_\alpha$ is defined using the topological pressure (rather than the Hentschel--Procaccia spectrum), and this family determines a ``predicted'' multifractal spectrum.  For ``non-mixed'' spectra, this expected spectrum is the Legendre transform of $T_0$.  For mixed spectra, the relationship is more subtle---details are given in Section~\ref{sec:mainresult}.  To keep our notation uniform, we denote the expected spectrum by $\SSS(\alpha)$ in both cases.

Now the multifractal formalism consists of showing that $\FFF(\alpha) = \SSS(\alpha)$.  In the Barreira--Saussol setting, this means establishing the equality $\dim_u K_\alpha = \SSS(\alpha)$, which is done by showing that $\SSS(\alpha)$ bounds $\dim_u K_\alpha$ both from above and from below.

As is typical of computations involving Hausdorff dimension and other dimensional quantities, the upper bound is easier to obtain than the lower bound.  Indeed, the inequality $\FFF(\alpha) \leq \SSS(\alpha)$ holds in almost complete generality; this is well-known for a number of spectra~\cite{oZ02,TB06,lO10}, and we show in Theorem~\ref{thm:universal} that it holds for $\FFF(\alpha)$ under minimal hypotheses.

Further hypotheses are necessary to obtain the lower bound $\FFF(\alpha) \geq \SSS(\alpha)$.  One may begin with the observation that on the interior of its domain (modulo a small technical condition), the predicted spectrum $\SSS(\alpha)$ can be written in terms of a conditional variational principle (Theorem~\ref{thm:Scvp}):
\[
\SSS(\alpha) = \sup \left\{ \frac{h_\mu(f)}{ \int u\,d\mu} \,\Big|\, \mu\in \Mf(X), 
\int (\ph - \alpha\psi) \,d\mu = 0, \int u\,d\mu > 0 \right\}.
\]
Given $\mu\in \Mf(X)$, we denote the set of $\mu$-generic points by $G_\mu := \{x\in X \mid \frac 1n S_n\phi(x) \to \int\phi\,d\mu \text{ for all } \phi\in C(X) \}$.  If $\int (\ph - \alpha \psi)\,d\mu=0$, then the ergodic theorem implies $G_\mu \subset K_\alpha$, and in particular, $\dim_u G_\mu \leq \dim_u K_\alpha$.  Writing 
\begin{equation}\label{eqn:Mfa}
\Mfa(X) := \left\{\mu\in \Mf(X) \,\Big|\, \int (\ph - \alpha\psi) \,d\mu = 0, \int u\,d\mu > 0 \right\},
\end{equation}
it follows that
\begin{equation}\label{eqn:cvps}
\sup_{\mu\in \Mfa(X)} (\dim_u G_\mu) \leq \FFF(\alpha) \leq
\sup_{\mu\in \Mfa(X)} \left( \frac{h_\mu(f)}{ \int u\,d\mu} \right).
\end{equation}
At this point there are two rather different approaches found in the literature to proving that equality holds in~\eqref{eqn:cvps}: the \emph{thermodynamic} approach and the \emph{saturation} approach.

\subsubsection{The thermodynamic approach}

For ergodic measures $\mu$, it was shown by Bowen~\cite{rB73} that $\htop G_\mu = h_\mu(f)$, and it is not hard to derive as a corollary that if $\int u\,d\mu > 0$, then 
\begin{equation}\label{eqn:saturated}
\dim_u G_\mu = \frac{h_\mu(f)}{\int u\,d\mu}.
\end{equation}
Thus equality holds in~\eqref{eqn:cvps} provided the second supremum does not change its value when we restrict our attention to ergodic measures.

An ergodic measure $\mu$ that achieves the second supremum in~\eqref{eqn:cvps} is called a \emph{full measure} for the spectrum $\FFF$ at the value $\alpha$.  In its simplest incarnation, the thermodynamic approach produces a full measure as an ergodic equilibrium state $\mu$ for a certain potential function $\phi$.  The fact that $\mu$ has the appropriate $u$-dimension and integrates $\ph-\alpha\psi$ to $0$ is guaranteed by using Bowen's equation (for $u$-dimension) and Ruelle's formula (for the derivative of pressure), under the assumption that the pressure function satisfies appropriate differentiability criteria at $\phi$. This is the approach followed in~\cite{BS01} and is also our primary strategy here.

For some spectra, there are values of $\alpha$ for which no full measure exists; this is often the case when a phase transition occurs (see Section~\ref{sec:MP}).  Nevertheless, it is often still possible to find a family of ergodic measures whose $u$-dimensions approach the second supremum in~\eqref{eqn:cvps}.  Corollary~\ref{cor:phase} gives conditions on the thermodynamics of a system that allow this to be done.

In addition to the equality between $\FFF(\alpha)$ and $\SSS(\alpha)$, the thermodynamic approach also establishes the stronger conditional variational principle
\begin{equation}\label{eqn:cvp3}
\dim_u K_\alpha = \sup \{ \dim_u \mu \mid \mu \in \Mfe(K_\alpha) \},
\end{equation}
where $\Mfe(X) := \{ \mu \in \Mf(X) \mid \mu \text{ is ergodic} \}$ and $\Mfe(K_\alpha) = \{\mu\in \Mfe(X) \mid \mu(K_\alpha) = 1\}$.

\subsubsection{The saturation approach}

One can also approach~\eqref{eqn:cvps} by extending~\eqref{eqn:saturated} to non-ergodic measures.  In general, $G_\mu$ need not even be non-empty when $\mu$ is not ergodic, and so extending~\eqref{eqn:saturated} requires further hypotheses on the map $f$.  As always, the difficult step is proving the lower bound $\dim_u G_\mu \geq h_\mu(f)/\int u\,d\mu$; the other bound is easy.

Most applications of the saturation approach in the literature assume that $f$ satisfies some version of the specification property and use this to show that $f$ is \emph{saturated}---that is, that $\htop G_\mu = h_\mu(f)$ for every invariant measure $\mu$~\cite{TV03,PS07}. 
The bound $\htop G_\mu \geq h_\mu(f)$ is obtained by using the specification property (or a variant) to build a measure $\nu$ supported on $G_\mu$ such that $\htop \nu := \inf \{\htop Z \mid\nu(Z)=1\} = h_\mu(f)$.  Because $G_\mu$ does not support any invariant measures, the measure $\nu$ is not $f$-invariant, and hence cannot be obtained via the usual thermodynamic methods, which only yield invariant measures as equilibrium states.

\subsubsection{Comparison of approaches and open problems}

The thermodynamic approach has the advantage that it does not require any knowledge of the system beyond its thermodynamic properties, and so can be applied to systems without the specification property or any other property known to imply saturation.

Furthermore, the thermodynamic approach establishes the strengthened conditional variational principle~\eqref{eqn:cvp3}, which the saturation approach does not prove, and which does not follow from having equality in~\eqref{eqn:cvps}.  This suggests the following problem.

\begin{problem}
Does there exist a dynamical system $f\colon X\to X$ and functions $\ph,\psi,u$ such that the associated multifractal spectrum $\FFF(\alpha)$ satisfies the multifractal formalism, but the conditional variational principle~\eqref{eqn:cvp3} fails and there exists $\alpha$ such that
\[
\sup \{ \dim_u \mu \mid \mu \in \Mfe(K_\alpha)\} < \dim_u K_\alpha = \sup \{ \dim_u \mu \mid \mu \in \Mfa(X) \}?
\]
\end{problem}

The saturation approach has the advantage that it works equally well in the presence of phase transitions---points of non-differentiability of the pressure function---whereas the thermodynamic approach requires careful application in these cases.  Indeed, if we restrict our attention to the case $u=1$, so $\dim_u = \htop$, then any saturated system satisfies the multifractal formalism for every pair of potentials $\ph,\psi \in C(X)$ with $\psi>0$.

\begin{problem}
Is saturation a necessary property for the full multifractal formalism to hold?  That is, suppose $f$ is not saturated; is it necessarily the case that there exist $\ph,\psi,u$ for which the formalism does not hold?
\end{problem}

Finally, the saturation approach can be used to establish detailed multifractal results about the irregular set (see Section~\ref{sec:irreg}), which cannot be studied using the usual thermodynamic approach, as it does not support any invariant measures.

\subsection{Outline of the paper}

In Section~\ref{sec:pressure}, we define the notions of topological pressure and $u$-dimension, which will be foundational to the rest of the paper.  Section~\ref{sec:mainresult} introduces the multifractal spectrum $\FFF(\alpha)$, the predicted spectrum $\SSS(\alpha)$, and gives our main results relating the two.  Section~\ref{sec:Scvp} shows that the predicted spectrum $\SSS(\alpha)$ can also be given in terms of a conditional variational principle.  Sections~\ref{sec:concave}--\ref{sec:dimspec} examine various spectra that can be obtained as specific instances of the main results.

Section~\ref{sec:survey} surveys some of the history of multifractal analysis and describes various results in the literature from both the thermodynamic and saturation approaches.  Section~\ref{sec:examples} gives applications of our main results to various classes of systems, including some non-uniformly hyperbolic ones.  Section~\ref{sec:pfs} contains all the proofs, and Section~\ref{sec:dim} contains background material on Carath\'eodory dimension characteristics from~\cite{yP98}.

\emph{Acknowledgments.}  Many of the results in this paper had their genesis in my Ph.D.\ dissertation, written under the supervision of Yakov Pesin.  I am grateful to him for his insightful guidance, unwavering support, and constant encouragement throughout.

\section{Definitions and results}\label{sec:results}

\subsection{Definition of pressure and $u$-dimension}\label{sec:pressure}

Before formulating our results precisely, we recall the definitions of topological pressure and $u$-dimension.  

\begin{definition}
Let $X$ be a compact metric space and $f\colon X\to X$ be continuous.  Fix $\phi\in C(X)$ and $Z\subset X$; denote by $S_n \phi(x)$ the Birkhoff sum $\phi(x) + \phi(f(x)) + \cdots + \phi(f^{n-1}(x))$.  For a given $\delta>0$ and $N\in\NN$, let $\PPP(Z,N,\delta)$ be the collection of countable sets $\{ (x_i,n_i) \} \subset Z\times \{N,N+1,\dots\}$ such that $Z \subset \bigcup_i B(x_i,n_i,\delta)$.  For each $s\in\RR$, consider the set functions
\begin{equation}\label{eqn:mZa}
\begin{aligned}
m_P(Z,s,\phi,N,\delta)&= \inf_{\PPP(Z,N,\delta)} \sum_{(x_i,n_i)} 
\exp\left(-n_i s + S_{n_i} \phi(x_i) \right), \\
m_P(Z,s,\phi,\delta)&=\lim_{N\to\infty} m_P(Z,s,\phi,N,\delta).
\end{aligned}
\end{equation}
This function is non-increasing in $s$, and takes values $\infty$ and $0$ at all but at most one value of $s$.  Denoting the critical value of $s$ by
\begin{align*}
P_Z(\phi,\delta) &= \inf \{s\in\RR \mid m_P(Z,s,\phi,\delta)=0\} \\
&= \sup \{s\in\RR \mid m_P(Z,s,\phi,\delta)=\infty\},
\end{align*}
we get $m_P(Z,s,\phi,\delta)=\infty$ when $s<P_Z(\phi,\delta)$, and $0$ when $s>P_Z(\phi,\delta)$.

The \emph{topological pressure} of $\ph$ on $Z$ is $P_Z(\phi) = \lim_{\delta\to 0} P_Z(\phi,\delta)$; the limit exists because given $\delta_1 < \delta_2$, we have $\PPP(Z,N,\delta_1) \subset \PPP(Z,N,\delta_2)$, and hence $m_P(Z,s,\phi,\delta_1) \geq m_P(Z,s,\phi,\delta_2)$, so $P_Z(\phi,\delta_1) \geq P_Z(\phi,\delta_2)$.  

In the particular case $\phi=0$, this definition yields the topological entropy $\htop(Z) = P_Z(0)$, as defined by Bowen for non-compact sets~\cite{rB73}.
\end{definition}

\begin{remark}
This definition differs formally from the one given by Pesin and Pitskel'~\cite{PP84} (see also~\cite{yP98}).  However, as shown in~\cite[Proposition 5.2]{vC10a}, the two definitions yield the same value.

When $Z=X$, we can also compute the pressure as follows.  For each $n\in \NN$ and $\delta>0$, let
\[
\ZZZ_n(\phi,\delta) = \inf \left\{ \sum_{x\in E} \exp\left(\sup_{y\in B(x,n,\delta)} S_n \phi(y)\right) \,\Big|\,
B(E,n,\delta) = X \right\},
\]
where we write $B(E,n,\delta) := \bigcup_{x\in E} B(x,n,\delta)$; then the pressure is given by
\[
P(\phi) = \lim_{\delta\to 0} \ulim_{n\to\infty} \frac 1n \log \ZZZ_n(\phi,\delta).
\]
\end{remark}

In the remainder of this paper, we will always assume that $\htop(X) < \infty$, since otherwise the pressure function always takes infinite values and does not contain the information we need.

We will consider multifractal spectra described by the following dimensional quantity, which was introduced by Barreira and Schmeling in~\cite{BS00}.

\begin{definition}
Let $X,f$ be as above and fix a continuous function $u\colon X\to \RR$ satisfying the following condition.
\begin{description}
\item[(P)]  For every $\delta>0$ there exist covers $E_N \in \PPP(X,N,\delta)$ such that $\lim_{N\to\infty} \inf_{(x,n)\in E_N} S_nu(x) = +\infty$.
\end{description}
In~\cite{BS00} it is assumed that $u>0$, which implies \PP.  However, there are examples where $u$ is not strictly positive but \PP\ is still satisfied (see Section~\ref{sec:MP}).

\begin{proposition}\label{prop:ugeq0}
If $u\in C(X)$ satisfies \PP, then $\int u\,d\mu \geq 0$ for every $\mu\in \Mf(X)$.
\end{proposition}

Given a set $Z\subset X$, consider the set functions
\begin{equation}\label{eqn:mua}
\begin{aligned}
m_u(Z,s,N,\delta)&= \inf_{\PPP(Z,N,\delta)} \sum_{(x_i,n_i)} e^{-s S_{n_i} u(x)}, \\
m_u(Z,s,\delta)&=\lim_{N\to\infty} m_u(Z,s,N,\delta).
\end{aligned}
\end{equation}
Once again, this function is non-increasing in $s$, and takes values $\infty$ and $0$ at all but at most one value of $s$.  Denoting the critical value of $s$ by
\[
\dim_u(Z,\delta) = \inf \{s\in\RR \mid m_u(Z,s,\delta)=0\},
\]
we get $m_u(Z,s,\delta)=\infty$ when $s<\dim_u(Z,\delta)$, and $0$ when $s>\dim_u(Z,\delta)$.

The \emph{$u$-dimension} of $Z$ is $\dim_u Z = \lim_{\delta\to 0} \dim_u(Z,\delta)$; the limit exists for the same reason as in the previous definition.  In the particular case $u=1$, this definition once again yields the topological entropy $\htop(Z)$.
\end{definition}

Topological pressure and $u$-dimension are both examples of Carath\'eodory dimension characteristics (see Section~\ref{sec:dim}).  Thus they both satisfy \emph{countable stability}: for any countable collection of arbitrary sets $Z_i \subset X$,
\begin{align}
\label{eqn:Pstable}
P_{\bigcup_i Z_i} (\phi) &= \sup_i P_{Z_i}(\phi), \\
\label{eqn:dimustable}
\dim_u \left( \bigcup_i Z_i \right) &= \sup_i \dim_u(Z_i).
\end{align}

An important tool for computing $u$-dimension is Bowen's equation.

\begin{proposition}\label{prop:bowen}
Let $X$ be compact, $f\colon X\to X$ be continuous, and $u\in C(X)$ satisfy \PP.  Suppose $Z\subset X$ has the property that $\llim_{n\to\infty} \frac 1n S_n u(x) > 0$ for all $x\in Z$.  Then
\begin{equation}\label{eqn:bowen}
\dim_uZ = \inf\{t\in \RR\mid P_Z(-tu) \leq 0\}.
\end{equation}
\end{proposition}

The most effective way to obtain lower bounds on $\dim_u Z$ is to study measures supported on $Z$.

\begin{definition}
The $u$-dimension of a probability measure $\mu \in \MMM(X)$ is
\begin{equation}\label{eqn:dimumu}
\dim_u \mu := \inf \{ \dim_u Z \mid \mu(Z) = 1 \}.
\end{equation}
\end{definition}

\begin{proposition}\label{prop:dimumu}
Given an ergodic $f$-invariant probability measure $\mu\in \Mfe(X)$ with $\int u\,d\mu > 0$, the $u$-dimension of $\mu$ is given by
\begin{equation}\label{eqn:dimumu2}
\dim_u \mu = \frac{h_\mu(f)}{\int_X u\,d\mu}.
\end{equation}
\end{proposition}

\begin{remark}
Note that~\eqref{eqn:dimumu2} fails for non-ergodic measures; indeed, the functions $\mu\mapsto h_\mu(f)$ and $\mu\mapsto \int u\,d\mu$ are both affine, while it is an easy consequence of~\eqref{eqn:dimumu} that $\dim_u(t\mu_1 + (1-t)\mu_2) = \max(\dim_u\mu_1, \dim_u\mu_2)$ for $0<t<1$.
\end{remark}

We emphasise once again that by the appropriate choice of the function $u$, the $u$-dimension $\dim_u$ can be made to equal either the topological entropy or the Hausdorff dimension (provided $f$ is conformal), and so all our results involving $\dim_u$ apply to both of these quantities.

\subsection{Mixed multifractal spectra and main results}\label{sec:mainresult}


Now we fix functions $\ph,\psi,u\in C(X)$ such that $u$ satisfies \PP.  Given $\alpha\in \RR$, consider the level set
\begin{equation}\label{eqn:level}
K_\alpha = K_\alpha(\ph,\psi) = \left\{ x\in X \,\Big|\, \lim_{n\to\infty} \frac{S_n \ph(x)}{S_n\psi(x)} = \alpha \right\}.
\end{equation}
The \emph{multifractal spectrum} associated to the triple $(\ph,\psi,u)$ is the function $\FFF\colon \RR\to \RR$ given by
\[
\FFF(\alpha) = \FFF_u^{\ph,\psi}(\alpha) := \dim_u K_\alpha(\ph,\psi),
\]
where we adopt the convention that $\dim_u \emptyset = -\infty$, and so $\FFF(\alpha)=-\infty$ for those values of $\alpha$ corresponding to empty level sets.

Consider the set of points
\begin{equation}\label{eqn:hatX}
\hat{X} := \left\{ x\in X \,\Big|\, \llim_{n\to\infty} S_n\psi(x) > 0 \text{ and } \llim_{n\to\infty} \frac 1n\log S_n u(x) > 0 \right\},
\end{equation}
and the level sets $\hat K_\alpha := K_\alpha \cap \hat X$.  We will obtain estimates on the spectrum
\[
\hat\FFF(\alpha) = \dim_u \hat K_\alpha(\ph,\psi).
\]
In our applications, we will find that $\hat\FFF = \FFF$ on the domain of interest, and thus our estimates apply directly to the original spectrum.

For every $\alpha\in \RR$, define a function $R_\alpha\colon \RR^2 \to \RR$ by 
\begin{equation}\label{eqn:R}
R_\alpha(q,t) := P(q(\ph - \alpha\psi) - tu),
\end{equation}
and define $T_\alpha\colon \RR\to \RR$ implicitly by
\begin{equation}\label{eqn:T}
T_\alpha(q) := \inf \{t \in \RR \mid R_\alpha(q,t)\leq 0\}.
\end{equation}
The \emph{thermodynamically predicted spectrum} $\SSS\colon \RR\to \RR$ is
\begin{equation}\label{eqn:S}
\SSS(\alpha) := \inf \{ T_\alpha(q) \mid q\in \RR \}.
\end{equation}
The predicted domain of the spectrum is
\begin{equation}\label{eqn:I}
I := \{ \alpha \in \RR \mid \SSS(\alpha) > -\infty \}.
\end{equation}
\begin{proposition}\label{prop:Sgeq0}
$\SSS(\alpha)\geq 0$ for every $\alpha\in I$.
\end{proposition}

In Section~\ref{sec:Scvp}, we will relate the predicted spectrum to a conditional variational principle, which is the form many results in the literature take.  In Section~\ref{sec:concave}, we will examine the particular case when $u=\psi$, which covers many of the most important spectra, and where $\SSS$ can be obtained as a Legendre transform.

Our first main result holds in great generality, and establishes one half of the multifractal formalism.

\begin{thma}\label{thm:universal}
Fix $\ph,\psi,u\in C(X)$ such that $u$ satisfies \PP.  Then 
\begin{equation}\label{eqn:FS}
\hat\FFF(\alpha) \leq \SSS(\alpha)
\end{equation}
for all $\alpha\in \RR$.  In particular, $\hat K_\alpha = \emptyset$ for all $\alpha\notin I$.
\end{thma}

Once again, we emphasise that in all our examples we have $\tilde\FFF=\FFF$, so that Theorem~\ref{thm:universal} applies directly to the spectrum of interest.  We also point out that~\eqref{eqn:FS} holds for $\alpha\in \di I$, but not for $\alpha=\infty$.

\begin{remark}
The infimum in~\eqref{eqn:S} may not be achieved by any finite value of $q$; Theorem~\ref{thm:Scvp} gives conditions under which it is achieved.  Even if it is achieved, however, the minimising value of $q$ may not be unique.
\end{remark}

A measure $\nu$ is a \emph{full measure} for the spectrum $\hat\FFF$ at $\alpha$ if $\nu$ is ergodic and invariant, $\nu(\hat K_\alpha) = 1$, and $\dim_u \nu = \SSS(\alpha)$.

\begin{thma}\label{thm:full-measure}
Let $\ph,\psi,u$ be as in Theorem~\ref{thm:universal}.  Fix $\alpha\in\RR$, and suppose that $q=q(\alpha)\in \RR$ realises the infimum in~\eqref{eqn:S}, so that $\SSS(\alpha) = T_\alpha(q)$.  Suppose also that 
\begin{enumerate}
\item the potential $q(\alpha)(\ph - \alpha\psi)-\SSS(\alpha)u$ has an equilibrium state $\nu_\alpha$ such that $\int \psi\,d\nu_\alpha>0$ and $\int u\,d\nu_\alpha>0$, and
\item the function $R_\alpha$ is differentiable 
at $(q(\alpha),\SSS(\alpha))$.
\end{enumerate}
Then $\nu_\alpha$ is a full measure for $\hat\FFF$ at $\alpha$; in particular, $\hat\FFF(\alpha)=\SSS(\alpha)$.
\end{thma}

Note that if $\hat\FFF(\alpha) = \FFF(\alpha)$ and Theorem~\ref{thm:full-measure} applies, then $\nu_\alpha$ is full for $\FFF$ as well, and $\FFF(\alpha)=\SSS(\alpha)$, so the multifractal formalism holds at $\alpha$.  This extends results of Barreira and Saussol~\cite{BS01}, which apply when the entropy map is upper semi-continuous, $\psi$ and $u$ are uniformly positive, and $R_\alpha$ is differentiable for all values of $\alpha,q,t$.

In Section~\ref{sec:MP}, we consider systems and spectra for which certain values of $\alpha$ do not have full measures, but the formalism still holds and $\FFF \equiv \SSS$.  For such systems, we need the following easy corollary of Theorem~\ref{thm:full-measure}.

\begin{corollary}\label{cor:phase}
Let $\ph,\psi,u$ be as in Theorem~\ref{thm:universal}.  Fix $\alpha\in\RR$, and suppose that there exist compact invariant sets $X_n\subset X$ such that the hypotheses of Theorem~\ref{thm:full-measure} hold at $\alpha$ for the restriction of the system to $X_n$.  Then equality holds in~\eqref{eqn:muFS}; in particular, $\hat\FFF(\alpha) = \SSS(\alpha)$.
\end{corollary}

\subsection{Conditional variational principle}\label{sec:Scvp}

As mentioned earlier, the predicted spectrum can be given in terms of a conditional variational principle for nearly every value of $\alpha$.

Given $\alpha\in \RR$, let $\KKK(\alpha) := \{ \int (\ph - \alpha\psi) \,d\mu \mid \mu\in\Mf(X) \}$, and let $\hat I := \{\alpha \mid 0\in \inter \KKK(\alpha) \}$.  Let $\Mfa(X)$ be as in~\eqref{eqn:Mfa}.

\begin{thma}\label{thm:Scvp}
Let $\ph,\psi,u$ be as in Theorem~\ref{thm:universal}.  Then for every $\alpha\in \hat I$, the infimum in~\eqref{eqn:S} is achieved at a finite value of $q$, and the predicted spectrum $\SSS(\alpha)$ satisfies
\begin{equation}\label{eqn:Scvp}
\SSS(\alpha) = \sup_{\mu\in \Mfa(X)} \left( \frac{h_\mu(f)}{ \int u\,d\mu} \right).
\end{equation}
\end{thma}

\begin{remark}
The equivalence between the two characterisations of $\SSS(\alpha)$ is well-known in many cases where the formalism holds (see~\cite{BS01}, for example).  Theorem~\ref{thm:Scvp} establishes this equivalence even when the formalism does not hold.
\end{remark}

\begin{proposition}\label{prop:Ihatint}
If there exists $\gamma>0$ such that $\int \psi\,d\mu \geq \gamma h_\mu(f)$ for all $\mu\in \Mf(X)$, then $\hat I = \inter I$.
\end{proposition}

\begin{remark}
The hypothesis of Proposition~\ref{prop:Ihatint} holds if $\psi$ is uniformly positive and $\htop (f) < \infty$; thus in this case Theorem~\ref{thm:Scvp} applies to all $\alpha\in \inter I$.
\end{remark}

The measures in the conditional variational principle~\eqref{eqn:Scvp} need not be supported on $K_\alpha$ unless they are ergodic.  Thus in general, we have the following extended version of~\eqref{eqn:FS}:
\begin{equation}\label{eqn:muFS}
\sup \left\{ \dim_u \mu \mid \mu \in \Mfe(K_\alpha) \right\} \leq \hat\FFF(\alpha) \leq \SSS(\alpha).
\end{equation}
It is of interest to know when all three quantities in~\eqref{eqn:muFS} are equal.  By Theorem~\ref{thm:full-measure}, this is the case when $R_\alpha$ is differentiable at $(q(\alpha),\SSS(\alpha))$, provided an equilibrium state $\nu$ exists with $\int \psi\,d\nu>0$ and $\int u\,d\nu>0$; furthermore, in this case the supremum is achieved.  Corollary~\ref{cor:phase} gives equality in~\eqref{eqn:muFS} under weaker conditions, but in this case the supremum may not be achieved.

\subsection{Concavity when $u=\psi$}\label{sec:concave}

There are several important examples of multifractal spectra for which $u=\psi$; these are the \emph{non-mixed} spectra, which have been more thoroughly studied.  For these spectra, $\SSS$ takes a simpler form and Theorem~\ref{thm:universal} can be strengthened.

With $\psi = u$, we have $R_\alpha(q,t) = P(q\ph - (\alpha q + t) u)$, and so $T_\alpha$ is given by
\[
T_\alpha(q) = \inf \{t\in \RR \mid P(q\ph - (\alpha q + t) u) \leq 0 \}.
\]
Observe that $T_0\colon \RR\to \RR$ is defined by
\begin{equation}\label{eqn:T0}
T_0(q) = \inf \{ t\in \RR \mid P(q\ph - tu) \leq 0 \},
\end{equation}
and we have $T_\alpha(q) = T_0(q) - \alpha q$.  In particular,
\begin{equation}\label{eqn:Leg}
\SSS(\alpha) = \inf \{ T_0(q) - \alpha q \mid q\in \RR \}
\end{equation}
is the Legendre transform of $T_0$.  Recall that the relationship of $T_0$ and $\SSS$ as a Legendre transform pair can be seen geometrically as follows:  for every $\alpha$, the value of $\SSS(\alpha)$ is the $y$-intercept of the subtangent line to the graph of $T_0$ with slope $\alpha$.  (This line is unique by convexity, although it need not intersect the graph at a unique point.)

\begin{proposition}\label{prop:Qfinite}
Let $Q := \{q\in \RR \mid T_0(q) < \infty \}$.  Then one of the following three cases holds:
\begin{enumerate}
\item $Q = \RR$;
\item $Q = (q_0,\infty)$ or $[q_0,\infty)$ for some $q_0\geq 0$;
\item $Q = (-\infty,q_0)$ or $(-\infty,q_0]$ for some $q_0\leq 0$.
\end{enumerate}
\end{proposition}

By convexity of $T_0$, the following limits exist for every $q\in \inter Q$:
\begin{align*}
D^-T_0(q) &:= \lim_{\eps\to 0^+} \frac 1\eps (T_0(q) - T_0(q-\eps)), \\
D^+T_0(q) &:= \lim_{\eps\to 0^+} \frac 1\eps (T_0(q+\eps) - T_0(q)).
\end{align*}
Writing $J(q) := [D^-T_0(q), D^+T_0(q)]$, we see that $\inter J(q) = \emptyset$ if and only if $T_0$ is differentiable at $q$.  Furthermore, because the sets $\inter J(q)$ are disjoint for different values of $q$ (a consequence of convexity), $\inter J(q)$ is non-empty for at most countably many values of $q$.  Observe that $\SSS$ is an affine function of $\alpha$ on each interval $J(q)$, and conversely, $\SSS$ is strictly concave on every interval that is disjoint from every $\inter J(q)$.

Now for every $q$, the graph of $T_0$ admits subtangents at $q$ with every slope in the interval $J(q)$.  Denote by $\tilde I$ the set of values of $\alpha$ corresponding to subtangents where $T_0$ is stably finite: $\tilde I := \bigcup_{q\in \inter Q} J(q)$.  Note that if $Q$ is open, then $\tilde I = I$.  However, when $Q=[q_0,\infty)$ or $Q=(-\infty,q_0]$, this need not be the case.

\begin{thma}\label{thm:concave-hull}
Let $\ph,\psi,u$ be as in Theorem~\ref{thm:universal}, and suppose that $u=\psi$.  Then $\SSS|_{\tilde I}$ is the concave hull of $\hat\FFF|_{\tilde I}$:
\[
\SSS(\alpha) = \inf \{ \tilde\SSS(\alpha) \mid \tilde\SSS \geq \hat\FFF \text{ is concave on $\tilde I$} \} \text{ for all } \alpha \in \tilde I.
\]
\end{thma}

\begin{remark}
When $u\neq \psi$, the spectrum $\SSS$ is not necessarily concave; examples where concavity fails are given in~\cite{BS01,IK09}.
\end{remark}

For the equivalence between the Legendre transform description of $\SSS$ and the conditional variational principle, we observe that in many important cases, we can apply Proposition~\ref{prop:Ihatint} to show that Theorem~\ref{thm:Scvp} applies on all of $\inter I$.

\begin{proposition}\label{prop:Ihatint2}
Let $X$ be a compact metric space, let $f\colon X\to X$ be continuous, and fix $\ph,\psi,u\in C(X)$ such that $\psi \equiv u$ satisfies \PP.  Suppose that every ergodic measure $\mu$ with positive entropy has $\int u\,d\mu>0$, and suppose further that $\dim_u X < \infty$.  Then we have $\int u\,d\mu \geq h_\mu(f) / \dim_u X$ for all $\mu\in \Mf(X)$, whence Proposition~\ref{prop:Ihatint} and Theorem~\ref{thm:Scvp} give~\eqref{eqn:Scvp} for all $\alpha\in \inter I$.
\end{proposition}

\subsection{Birkhoff spectrum and coarse spectra}\label{sec:coarse}

We temporarily restrict our attention to the case where $\psi = u \equiv 1$.  In this case the predicted spectrum takes on a particularly simple form:  indeed, we see that
\[
T_0(q) = \inf \{t\in \RR \mid P(q\ph - t) \leq 0 \} = \inf \{t\in \RR \mid P(q\ph) - t \leq 0 \}
= P(q\ph),
\]
and so $\SSS$ is the Legendre transform of the convex function $q\mapsto P(q\ph)$.

For $\psi = u \equiv 1$, the level sets are
\begin{equation}\label{eqn:level2}
K_\alpha^\ph = \left\{ x\in X \,\Big|\, \lim_{n\to\infty} \frac 1n S_n \ph(x) = \alpha \right\}
\end{equation}
and the multifractal spectrum associated to $(\ph,\psi,u)$ is the entropy spectrum for Birkhoff averages
\[
\FFF(\alpha) = \BBB_E(\alpha) = \htop(K_\alpha^\ph).
\]
We will refer to this as the \emph{fine spectrum} to distinguish it from the \emph{coarse spectrum}, which we now define.

Given a neighbourhood $U\subset \RR$ of $\alpha$, consider the approximate level set
\begin{equation}\label{eqn:approx}
G_n(U) := \left\{ x\in X \,\Big|\, \frac 1n S_n \ph(x) \in U \right\},
\end{equation}
and let
\begin{equation}\label{eqn:Nn}
\Lambda_n^\delta(U) := \min \{ \# E \mid B(E,n,\delta) \supset G_n(U) \},
\end{equation}
where $B(E,n,\delta) := \bigcup_{x\in E} B(x,n,\delta)$.  The coarse multifractal spectra are the functions $\lF,\uF \colon \RR \to \RR$ given by
\begin{align*}
\lF(\alpha) &:= \inf_{U\ni \alpha} \lim_{\delta\to 0} \llim_{n\to\infty} \frac 1n \log \Lambda_n^\delta(U), \\
\uF(\alpha) &:= \inf_{U\ni \alpha} \lim_{\delta\to 0} \ulim_{n\to\infty} \frac 1n \log \Lambda_n^\delta(U).
\end{align*}

In this setting, we have the following stronger version of Theorem~\ref{thm:universal}.

\begin{thma}\label{thm:universal2}
Let $\ph\in C(X)$.  Then for all $\alpha$,
\begin{equation}\label{eqn:FFS}
\FFF(\alpha) \leq \lF(\alpha) \leq \uF(\alpha) \leq \SSS(\alpha).
\end{equation}
\end{thma}

For convenience and concreteness, we state a corollary that collects our main results in the setting where $\psi = u = 1$ and the entropy map is upper semi-continuous (which implies existence of equilibrium states).

\begin{corollary}[Entropy spectrum for Birkhoff averages]\label{cor:birkhoff}
Let $X$ be a compact metric space, $f\colon X\to X$ be continuous, and $\ph\in C(X)$.  Suppose that the entropy map $\Mf(X) \to \RR$ is upper semi-continuous.  Define $T_0(q) = P(q\ph)$ and let $J(q) = [D^-T_0(q),D^+T_0(q)]$.  Then the predicted spectrum is $\SSS(\alpha) = \inf \{ T_0(q) - q\alpha \mid q\in \RR \}$ and the predicted domain $I=\{\alpha\mid \SSS(\alpha)\geq 0\}$ is given by
\begin{equation}\label{eqn:Ibirk}
I = \left[\lim_{q\to-\infty} D^-T_0(q), \lim_{q\to\infty} D^+T_0(q)\right] = \left\{ \int \ph\,d\mu \,\Big|\, \mu\in \Mf(X) \right\}.
\end{equation}
For every $\alpha\notin I$ we have $K_\alpha^\ph=\emptyset$, while for every $\alpha\in \inter I \setminus \bigcup_{q\in \RR} \inter J(q)$, we have
\begin{equation}\label{eqn:BEcvp}
\begin{aligned}
\htop (K_\alpha^\ph) &= \lF(\alpha) = \uF(\alpha) \\
&= \inf \{ T_0(q) - q\alpha \mid q\in \RR \} \\
&= \sup \left\{ h_\mu(f) \,\Big|\, \mu\in \Mf(X), \int \ph\,d\mu = \alpha \right\} \\
&= \sup \{ h_\mu(f) \mid \mu \in \Mfe(K_\alpha^\ph)\}.
\end{aligned}
\end{equation}
Finally, if there exist compact invariant sets $X_n\subset X$ such that $P_{X_n}(q\ph) \to P_X(q\ph)$ for all $q$ and $q\mapsto P_{X_n}(q\ph)$ is differentiable on $\RR$ for all $n$, then~\eqref{eqn:BEcvp} holds for every $\alpha\in \inter I$.
\end{corollary}

\subsection{Conformal maps and Lyapunov spectra}\label{sec:Lyapunov}

The map $f\colon X\to X$ is \emph{conformal} with conformal derivative $f'(x)$ if for every $x\in X$ we have
\begin{equation}\label{eqn:conformal}
f'(x) = \lim_{y\to x} \frac{d(f(x),f(y))}{d(x,y)},
\end{equation}
where $f'\colon X\to [0,\infty]$ is continuous.  We restrict our attention to maps without critical points or singularities---that is, maps for which $0 < f'(x) <\infty$ for all $x\in X$.

When $f$ is conformal, the Lyapunov exponent of a point $x$ is $\lambda(x) = \lim \frac 1n S_n \log f'(x)$, provided the limit exists, and the level sets for Lyapunov exponents are
\[
K_\alpha^\lambda = \{ x\in X \mid \lambda(x) = \alpha \}.
\]
One important spectrum is the \emph{entropy spectrum for Lyapunov exponents} $\LLL_E$, which is the multifractal spectrum associated to $\ph = \log f'$ and $\psi = u = 1$.  This is nothing but the Birkhoff spectrum for the potential $\ph$, and so has already been discussed in Section~\ref{sec:coarse}.

One may also consider the \emph{dimension spectrum for Lyapunov exponents} $\LLL_D$, which is obtained by replacing $\htop K_\alpha$ with $\dim_H K_\alpha$.  This spectrum can also be obtained in the present framework by replacing $u=1$ with $u=\log f'$.  Indeed, when $u=\log f'$ is uniformly positive (that is, when $f$ is uniformly expanding), it was shown in~\cite{BS00} that $\dim_u Z = \dim_H Z$ for all $Z\subset X$, so the $u$-dimension in this case coincides with the Hausdorff dimension.

The relationship between $u$-dimension and Hausdorff dimension for non-uniformly expanding systems is slightly more delicate.  We consider sets all of whose points have positive lower Lyapunov exponent and satisfy the following \emph{tempered contraction} condition for $u=\log f'$:
\begin{equation}\label{eqn:tempered}
\inf_{0\leq k \leq n < \infty} \{ S_{n-k} u(f^k(x)) + n\eps \} > -\infty \text{ for every } \eps>0.
\end{equation}
In our examples, we will have $f'(x)\geq 1$ everywhere, which implies that every $x\in X$ satisfies~\eqref{eqn:tempered}.  For the time being, however, we allow $f'(x) < 1$.

Recall that our results apply to level sets in the set $\hat X$ given in~\eqref{eqn:hatX}.  With $\psi \equiv 1$, we have that $\hat{X} =  \{ x \in X \mid \llim_{n\to\infty} \frac 1n S_n u(x) > 0 \}$; consider the set $X^\star = \{x\in \hat X \mid \text{\eqref{eqn:tempered} holds} \}$.  It is shown in~\cite{vC10a} that for every $Z\subset X^\star$, we have
\[
\dim_H Z = \inf \{ t \geq 0 \mid P_Z(-tu) \leq 0 \}.
\]
In particular, using Proposition~\ref{prop:bowen}, we see that every such set has $\dim_u Z = \dim_H Z$.

It is shown in~\cite[Proposition~5.5]{vC10a} that if the Lyapunov exponent of $x$ exists---that is, if $\lim_{n\to\infty} \frac 1n S_n u(x)$ exists---then $x$ satisfies~\eqref{eqn:tempered}.  Thus if $Z\subset X$ is such that the Lyapunov exponent exists and is positive for every $x\in Z$, then $\dim_u Z = \dim_H Z$.  In particular, this shows that for \emph{any} conformal map (not necessarily uniformly expanding) without critical points or singularities, the mixed multifractal spectrum associated to $\ph = u = \log f'$ and $\psi = 1$ is equal to the dimension spectrum for Lyapunov exponents, at least for $0 < \alpha < \infty$.  (The case $\alpha=0$ is more subtle.)  This quickly leads to the relationship $\LLL_D(\alpha) = \frac 1\alpha \LLL_E(\alpha)$ for every conformal expanding map without critical points or singularities.

\begin{corollary}[Dimension spectrum for Lyapunov exponents]\label{cor:Lyapunov}
Let $X$ be a compact metric space, and let $f\colon X\to X$ be continuous and conformal without critical points or singularities, and suppose that $u=\log f'$ satisfies \PP.  Suppose also that the entropy map $\Mf(X) \to \RR$ is upper semi-continuous, and write $\lambda(\mu) = \int u\,d\mu = \int \log f'\,d\mu$ for the Lyapunov exponent of a measure $\mu$.

Write $T_0(q) := P(q\log f')$, and let $J(q) = [D^-T_0(q), D^+T_0(q)]$.  Then for all $\alpha>0$, the predicted spectrum is $\SSS(\alpha) = \frac 1\alpha \inf \{ T_0(q) - q\alpha \mid q\in \RR \}$ and the predicted domain $I=\{\alpha\mid \SSS(\alpha)\geq 0\}$ is given by
\begin{equation}\label{eqn:ILyap}
I = \left[\lim_{q\to-\infty} D^-T_0(q), \lim_{q\to\infty} D^+T_0(q)\right] = \{ \lambda(\mu) \mid \mu\in \Mf(X)\}.
\end{equation}
For every $\alpha\notin I$ we have $K_\alpha^\lambda=\emptyset$, while for every $\alpha\in \inter I \setminus \bigcup_{q\in \RR} \inter J(q)$, we have
\begin{equation}\label{eqn:LDcvp}
\begin{aligned}
\dim_H (K_\alpha^\lambda) &= \frac 1\alpha \inf \{ T_0(q) - q\alpha \mid q\in \RR \} \\
&= \sup \left\{ \frac{h_\mu(f)}{\lambda(\mu)} \,\Big|\, \mu\in \Mf(X), \lambda(\mu) = \alpha \right\} \\
&= \sup \left\{ \frac{h_\mu(f)}{\lambda(\mu)} \,\Big|\, \mu \in \Mfe(K_\alpha^\lambda)\right\}.
\end{aligned}
\end{equation}
Finally, if there exist compact invariant sets $X_n\subset X$ such that $P_{X_n}(q\ph) \to P_X(q\ph)$ for all $q$ and $q\mapsto P_{X_n}(q\ph)$ is differentiable on $\RR$ for all $n$, then~\eqref{eqn:LDcvp} holds for every $\alpha\in \inter I$.
\end{corollary}

It is worth mentioning that because $u\neq \psi$, the dimension spectrum for Lyapunov exponents is not a Legendre transform and need not be concave.  Examples where concavity fails have been given in~\cite{IK09}.

\subsection{Dimension spectrum for weak Gibbs measures}\label{sec:dimspec}

Let $\mu$ be a Borel probability measure on a separable metric space $X$.  The \emph{pointwise dimension} of $\mu$ at $x\in X$ is
\[
d_\mu(x) := \lim_{r\to 0} \frac{\log \mu(B(x,r))}{\log r},
\]
provided the limit exists.  Denoting the level sets by $K_\alpha^{d_\mu} := \{x\in X \mid d_\mu(x) = \alpha\}$, the \emph{dimension spectrum for pointwise dimensions} of $\mu$ is
\[
\DDD_D(\alpha) := \dim_H K_\alpha^{d_\mu}.
\]
This definition does not require a dynamical system in the background.  However, we will focus on the case where $\mu$ is a weak Gibbs measure for a potential $\phi\in C(X)$ and a conformal map $f\colon X\to X$---that is, where
\begin{equation}\label{eqn:wkgibbs}
\begin{aligned}
P(\phi) &= \lim_{\delta\to 0} \llim_{n\to\infty} \frac 1n (-\log \mu(B(x,n,\delta)) + S_n \phi(x)) \\
&= \lim_{\delta\to 0} \ulim_{n\to\infty} \frac 1n (-\log \mu(B(x,n,\delta)) + S_n \phi(x))
\end{aligned}
\end{equation}
for every $x\in X$.  (Note that $\mu$ is not required to be invariant.)  The dimension spectrum $\DDD_D$ is obtained in the present framework by taking $\psi = u = \log f'$ and $\ph=P(\phi)-\phi$.  Observe that from~\eqref{eqn:wkgibbs}, we have $\llim \frac 1n S_n \ph(x) \geq 0$ for all $x\in X$.\footnote{Our choice of $\ph$ differs from some other papers in the literature, which use $\ph = \phi-P(\phi)$.  We use the opposite sign so as to maintain consistency between the form of the Legendre transform in~\eqref{eqn:DDcvp} and that already given in~\eqref{eqn:BEcvp}.}

We consider non-uniformly expanding systems with $f'(x)\geq 1$ for all $x\in X$; thus $u(x) \geq 0$.  It is often the case that even though we may have $\llim \frac 1n S_n u(x) = 0$ for some $x\in X$, we nevertheless have
\begin{equation}\label{eqn:uphi}
\llim_{n\to\infty} \frac 1n S_n (u + \ph)(x) > 0
\end{equation}
for all $x\in X$.  (See Section~\ref{sec:MP}.)

\begin{corollary}[Dimension spectrum for pointwise dimensions]\label{cor:dimspec}
Let $X$ be a compact metric space, and let $f\colon X\to X$ be continuous and conformal with $1\leq f'(x) < \infty$ for all $x\in X$.  Suppose that the entropy map $\Mf(X) \to \RR$ is upper semi-continuous.  

Let $\mu$ be a weak Gibbs measure in the sense of~\eqref{eqn:wkgibbs} for $\phi\in C(X)$, and set $\ph = P(\phi) - \phi$ and $u=\log f'$.  Suppose that~\eqref{eqn:uphi} holds for all $x\in X$ and that $u$ satisfies \PP.  Furthermore, suppose that $\dim_u X<\infty$ and every measure $\nu\in \Mf(X)$ with $\lambda(\nu) := \int u\,d\nu = 0$ has $h_\nu(f) = 0$.  Let
\[
T_0(q) := \inf \{t\in \RR \mid P(q\ph - t\log a) \leq 0 \},
\]
and let $J(q) = [D^-T_0(q), D^+T_0(q)]$.  Then the predicted spectrum is $\SSS(\alpha) = \inf \{ T_0(q) - q\alpha \mid q\in \RR \}$ and the predicted domain $I=\{\alpha\mid \SSS(\alpha)\geq 0\}$ is given by
\begin{equation}\label{eqn:Idim}
I = \left[\lim_{q\to-\infty} D^-T_0(q), \lim_{q\to\infty} D^+T_0(q)\right].
\end{equation}
For every $\alpha\notin I$ we have $K_\alpha^{d_\mu}=\emptyset$, while for every $\alpha\in \inter I \setminus \bigcup_{q\in \RR} \inter J(q)$, we have
\begin{equation}\label{eqn:DDcvp}
\begin{aligned}
\dim_H (K_\alpha^{d_\mu}) &= \inf \{ T_0(q) - q\alpha \mid q\in \RR \} \\
&= \sup \left\{ \frac{h_\nu(f)}{\lambda(\nu)} \,\Big|\, \nu\in \Mf(X), \frac{\int \ph\,d\nu}{\lambda(\nu)} = \alpha \right\} \\
&= \sup \left\{ \frac{h_\nu(f)}{\lambda(\nu)} \,\Big|\, \nu \in \Mfe(K_\alpha^{d_\mu})\right\}.
\end{aligned}
\end{equation}
Finally, if there exist compact invariant sets $X_n\subset X$ such that $P_{X_n}(q\ph) \to P_X(q\ph)$ for all $q$ and $(q,t)\mapsto P_{X_n}(q\ph - tu)$ is differentiable on $\RR^2$ for all $n$, then~\eqref{eqn:DDcvp} holds for every $\alpha\in \inter I$.
\end{corollary}

\begin{remark}
In fact, Corollary~\ref{cor:dimspec} holds under slightly more relaxed assumptions on $f'$; namely, we can have $f'(x)<1$ provided $f$ has no critical points and~\eqref{eqn:tempered} holds for every $x\in K_\alpha^{\ph,u} \cup K_\alpha^{d_\mu}$ for which $\llim \frac 1n S_n u(x) > 0$.  Similarly, the condition that~\eqref{eqn:uphi} holds for all $x\in X$ can be relaxed to the condition that it holds for every $x\in K_\alpha^{\ph,u} \cup K_\alpha^{d_\mu}$ with $0\leq \alpha<\infty$.
\end{remark}

\section{History and survey of known results}\label{sec:survey}

In this section we describe some of the settings in which the multifractal formalism has been established to date.  Surveys of the literature can also be found in~\cite{yP98,lB08}, among others.

\subsection{Results for $\DDD_D$}\label{sec:DD}

The dimension spectrum for pointwise dimensions has received the most attention of any multifractal spectrum.  It was introduced in~\cite{HJKPS86}, although~\cite{BPPV84} already considered similar notions in the study of turbulence.

In order to establish the multifractal formalism for the dimension spectrum of a measure $\mu$, some information about the structure of $\mu$ is required.  In particular, $\mu$ is usually assumed to have some sort of Gibbs property that relates its scaling behaviour to the Birkhoff averages of a potential function $\ph$ under a map $f$.  The easiest case is when $f$ is uniformly expanding and conformal---the latter hypothesis is automatically fulfilled if we study one-dimensional differentiable maps (real or complex). 

We first mention the relevant papers for the uniformly expanding case.  Maps of the interval satisfying a Markov property were studied in~\cite{CLP87}, which dealt with the dimension spectrum for Markov measures, and also in~\cite{dR89}, which dealt with more general Gibbs measures.  In the complex plane, rational maps that are uniformly expanding on their Julia set were studied in~\cite{aL89}, which established the formalism for the measure of maximal entropy.  Finally, in~\cite{PW97}, Pesin and Weiss established the formalism for arbitrary conformal repellers and Gibbs measures corresponding to H\"older continuous potentials.  They use the thermodynamic approach to produce a family of full measures for the dimension spectrum; many of their results follow from Corollary~\ref{cor:dimspec} of the present paper (see Section~\ref{sec:PW}), and also from the results in~\cite{BS01}, since for uniformly expanding conformal maps it was shown in that paper that $\dim_u = \dim_H$ when $u=\log f'$.

The dimension spectrum for conformal measures on piecewise monotonic interval maps (with no assumption of any Markov property) was studied by Hofbauer~\cite{fH95}, using thermodynamic techniques based on approximating $f$ with a Markov map and $\ph,u$ by locally constant functions.  Despite the thermodynamic nature of this approach, the approximation techniques permit stronger results that apply even when the pressure function is non-differentiable.

The dimension spectrum for Gibbs measures has also been studied for maps with both expanding and contracting directions, provided the dynamics on each of those directions are conformal (which is the case if the expanding and contracting distributions are each one-dimensional).  The geodesic flow on a surface of constant negative curvature was studied in~\cite{oR95}, while Axiom A surface diffeomorphisms were studied in~\cite{dS94a,dS94b,PW97b}, and conformal Axiom A flows were studied in~\cite{PS01}.

Prior to this, the dimension spectrum for measures arising as a result of a geometric construction (rather than as invariant measures for a dynamical system) was studied in~\cite{BMP92, CM92, EM92}, where the thermodynamically defined function $T_0$ is replaced by a function related to the Hentschel--Procaccia dimension spectrum for $\mu$.  This line of inquiry has been vigorously continued---see~\cite{kF94,LN98,LN99,HMU02,MU03,bT06}, among many others.  In this paper, we will focus on measures arising as invariant measures for dynamical systems.

Within this context, there are three primary ways to weaken the hypotheses on the triple $(\mu,f,\ph)$.  In the first place, one can consider potentials $\ph$ that are not H\"older continuous, which often requires weakening the Gibbs property for $\mu$.  Secondly, one can replace the uniform hyperbolicity of $f$ with something weaker.  Thirdly (and most challengingly), one can remove the requirement that $f$ be conformal.

In the symbolic setting, the dimension spectrum for a class of measures called $g$-measures was studied in~\cite{eO99,eO00}, and weak Gibbs measures for continuous circle maps were studied in~\cite{FO03}.  These examples include cases where the pressure function is not differentiable, and the proofs rely on saturation techniques.

The easiest class of non-uniformly hyperbolic maps to consider are expanding Markov interval maps with parabolic fixed points, such as the Manneville--Pomeau map $x\mapsto x+x^{1+\beta} \pmod 1$, where $\beta>0$, which has a parabolic fixed point at $0$.  On the symbolic models for these systems, partial results for the dimension spectrum of equilibrium states of suitably regular measures were obtained in~\cite{PW99}.  Further results for the dimension spectrum of the measure of maximal entropy for the original system were given in~\cite{kN00}, and more general weak Gibbs measures were studied in~\cite{mY02}.  In the invertible regime, horseshoe maps with parabolic fixed points were considered in~\cite{BI10}.

All the results in the previous paragraph make use of inducing schemes and the thermodynamic approach;  see Section~\ref{sec:MP} for a description of similar results that can be obtained using Corollary~\ref{cor:dimspec} (without using inducing schemes).

For Markov interval maps with a parabolic fixed point, the dimension spectrum for weak Gibbs measures associated to continuous potentials (a broader class of potentials than those above) was studied in~\cite{JR09} using what are essentially saturation techniques.

One-dimensional maps with critical points are more difficult to study; in this setting, the presence of the critical point prevents a direct application of Corollary~\ref{cor:dimspec}, and so results beyond those given in this paper are needed.  Using inducing schemes and the corresponding thermodynamics, the dimension spectra of various conformal measures have been examined in~\cite{IT09, mT09}.

Non-conformal maps are another kettle of fish entirely, requiring new methods; in particular, the classical thermodynamic formalism must be replaced with the non-additive formalism introduced in~\cite{lB96}.  Dimension spectra for various non-conformal maps have been studied in~\cite{jK95,lO98,JR09b}.

Finally, we mention that the relationship given in Theorem~\ref{thm:concave-hull} between the multifractal spectrum and a globally defined function has been presented already in~\cite{TB06}, which shows that the Legendre transform of a certain free energy function gives the concave hull of $\DDD_D(\alpha)$ in a very general setting.  The free energy function in question is given not in terms of the topological pressure, but in terms of the measure $\mu$ being studied, and is related to the Hentschel--Procaccia spectrum.

\subsection{Results for $\LLL_D$}

The dimension spectrum for Lyapunov exponents was studied for one-dimensional full-branched Markov maps by Bohr and Rand~\cite{BR87}.  The formalism was extended to Axiom A surface diffeomorphisms in~\cite{PW97b}, to conformal repellers in~\cite{hW99}, and to conformal Axiom A flows in~\cite{PS01}, using thermodynamic techniques in each case.  See Section~\ref{sec:PW} for applications of Corollary~\ref{cor:Lyapunov} in some of these settings; because the maps in question are uniformly expanding, the results in~\cite{BS01} also apply here.

More recently, using similar thermodynamic and approximation techniques to those in~\cite{fH95}, Hofbauer studied piecewise monotonic interval maps in~\cite{fH10b}; note that the lack of any Markov structure makes this case significantly more difficult than that studied in~\cite{BR87}.

Markov maps of the interval with parabolic fixed points were studied in~\cite{PW99}, where results were obtained on part of the domain of the spectrum;  complete results in this setting were given in~\cite{kN00}.  These proofs use inducing schemes and the thermodynamic approach; see Section~\ref{sec:MP} for some applications of Corollary~\ref{cor:Lyapunov} in this non-uniformly hyperbolic context.

Interval maps with parabolic fixed points were studied again in~\cite{GR09}, and rational maps whose Julia sets contain parabolic points were considered in~\cite{GPR10,GPRR10}.  In addition to determining the Hausdorff dimension of level sets, these papers use the saturation approach to obtain detailed information about the Hausdorff dimension of the irregular sets $\{x \mid \llambda(x)=\alpha, \ulambda(x)=\beta \}$, as well as the set of points with lower Lyapunov exponent equal to $0$.  

As with dimension spectra, the fully non-uniform case of multimodal maps cannot be treated with the results in this paper.  The Lyapunov spectrum for multimodal maps has been studied by Iommi and Todd~\cite{IT09} using inducing schemes and thermodynamic results.

For non-conformal maps, life is again more difficult; some results on the \emph{entropy} spectrum for Lyapunov exponents are given in~\cite{BG06}.

\subsection{Results for $\BBB_E$}

The entropy spectrum for Birkhoff averages (or something equivalent to it) of arbitrary continuous potentials was studied using saturation techniques for shifts of finite type in~\cite{FF00,FFW01}.  This was extended to systems with specification in~\cite{TV03,FLP08}, and to the weaker $g$-almost product property in~\cite{PS07}.  For systems with specification, a saturation property for the pressure function was shown in~\cite{dT09}.

The thermodynamic approach was initially applied to systems with specification and sufficiently regular potential functions in~\cite{TV99}, and some results on subshifts of finite type are given in~\cite{PW01}.  In~\cite{BS01}, it was observed that differentiability of the pressure function $q\mapsto P(q\ph)$ is the key, leading to a version of Corollary~\ref{cor:birkhoff}.  A similar result for flows was shown in~\cite{BD04}.  

More recently, a version of Corollary~\ref{cor:birkhoff} for almost additive potential functions was proved in~\cite{FH10}; non-additive versions of the Birkhoff spectrum have also been studied in~\cite{BD09}, and a broad class of non-additive functions obtained from matrix-value cocycles are studied in~\cite{dF03,dF09,BT09}.


\subsection{Results for mixed spectra}

One of the simplest mixed spectra to study is the dimension spectrum for Birkhoff averages.  For uniformly expanding conformal maps, we have $\dim_u = \dim_H$ when $u=\log f'$, and so this spectrum was studied in~\cite{BS01}, using thermodynamic techniques and assuming that the pressure function $(q,t)\mapsto P(q\ph - t\log f')$ is differentiable.  The corresponding results for flows were proved in~\cite{BD04}.  

For non-uniformly expanding maps with parabolic fixed points, this spectrum was studied using saturation techniques in~\cite{JJOP08}.  Because of the use of saturation techniques rather than thermodynamic techniques, the results obtained apply to all continuous potential functions, and not just those for which the pressure function is differentiable.

See Section~\ref{sec:MP} for results in the non-uniformly expanding case using the thermodynamic approach and applying Corollary~\ref{cor:phase} (and hence establishing the stronger variational principle).  As before, these results do not yet extend to the dimension spectrum for Birkhoff averages on multimodal maps, which is considered in~\cite{yC08b}.

Various non-conformal systems and their associated dimension spectra for Birkhoff averages have been studied in~\cite{BR07,JS07,hR10,hR10b}.

More general mixed spectra, such as the entropy spectrum for pointwise dimensions and the dimension spectrum for local entropies, were mentioned in the general framework of~\cite{BPS97}, but no results were proven.  The $\dim_u$-spectra studied in this paper were introduced in~\cite{BS01}, where the formalism was proven under assumptions of differentiability of the pressure function.

\subsection{Irregular set}\label{sec:irreg}

For the local asymptotic quantities used in defining a multifractal spectrum, it is usually the case that the irregular set of points on which the lower and upper limits disagree is invisible to invariant measures.  Thus thermodynamic methods do not suffice, and one needs saturation techniques to show that the irregular set has full dimension, which is generally the case provided it is non-empty.

For mixing shifts of finite type, the irregular set for Birkhoff averages was shown to have full entropy in~\cite{FFW01}.  A more general result for $u$-dimensions and shifts with specification was shown in~\cite{BS00,BS01}, and in~\cite{dT11}, it is shown that the irregular set carries full topological pressure for systems with specification.

For uniformly expanding conformal maps with a Markov structure, the Hausdorff and packing dimensions of various sorts of irregular sets, including more refined sets defined as those points $x$ for which the quantities $\phi_n(x)$ have a specified accumulation set, are studied in the series of papers~\cite{lO03,OW07,lO06,lO08}.  For systems with specification, analogous results for topological entropy are given in~\cite{EKL05}.

For maps with parabolic fixed points, the irregular set for Lyapunov exponents was studied in~\cite{PW99} and more recently in~\cite{GR09,GPR10,GPRR10}.  The latter set of papers contains detailed information about not just the irregular set itself, but the set of points whose sequence of approximate Lyapunov exponents has a specified lower and upper limit.

For multimodal interval maps, the Hausdorff dimension of the irregular set for Birkhoff averages was studied in~\cite{yC08b}, and the irregular set for pointwise dimensions was studied in~\cite{mT09}.


Another sort of irregular set, comprising points for which the forwards and backwards averages disagree, was studied in~\cite{mS91} for hyperbolic surface diffeomorphisms.

\subsection{Coarse spectra}\label{sec:coarse2}

Coarse multifractal spectra have received rather less attention than their fine counterparts, and what attention they have received has been focused predominantly on the coarse spectrum for pointwise dimensions.  This spectrum was studied in various settings in~\cite{kF03,RM97,RM98,lO05,TB06}.  An analogue of Theorem~\ref{thm:concave-hull} has been given in~\cite{lO10}, showing that the Legendre transform of the coarse spectrum for dimensions is the R\'enyi spectrum.

\subsection{Higher-order spectra}\label{sec:multidim}

A natural and important generalisation of the multifractal theory surveyed in this paper is to allow the local quantities $\phi_n(x)$ to take values in a more general space than $\RR$.  In its simplest form, replacing $\RR$ with $\RR^2$, this amounts to considering the spectrum $\FFF\colon \RR^2\to \RR$ defined by
\[
\FFF(\alpha,\beta) = \dim_u (K_\alpha\cap K_\beta).
\]
The results from the one-dimensional case do not immediately apply to these higher-order spectra, which exhibit various new phenomena: in particular, the domain of the spectrum need not be convex, and its interior may be disconnected~\cite{BSS02}.

Nevertheless, many of the same techniques can be used in this higher-order setting.  The thermodynamic approach was applied in \cite{BSS02} to show that differentiability of the appropriate cross-section of the pressure function still leads to a conditional variational principle (see also \cite{lB08}).  

Using saturation techniques, higher-order spectra for subshifts of finite type were studied in~\cite{FF00,FFW01}, and also~\cite{lO03,OW07,lO06,lO08}.  (The latter series of papers also studies more general local quantities.)  For systems with specification, multifractal spectra for quantities taking values in $\RR^d$ were studied in~\cite{TV03}, and Banach-space valued quantities in~\cite{FLP08}.

\subsection{Other multifractal results}

There are a great many papers in the literature addressing other questions of a multifractal nature that are closely related to the ones mentioned so far.  We mention just a few of these.

In~\cite{lO04c, JP07}, the level sets $K_\alpha$ are further partitioned into sets $K_{\alpha,\gamma}$, where $\gamma$ controls the \emph{rate} of convergence $\phi_n(x)\to \alpha$, and the sizes of these sets are examined.

In~\cite{gI05,gI10,IJ10}, various multifractal spectra are studied in the context of a topological Markov chain with \emph{infinitely} many symbols; this is related to iterated function systems with countably many branches~\cite{KU07}.

In~\cite{BQ10}, the asymptotic behaviour captured by the level sets $K_\alpha$ is localised: fixing $\xi\in C(X)$, one studies the multifractal spectrum defined with the level sets $K_\alpha = \{x\in X \mid \phi_n(x) \to \xi(x) \}$.

In~\cite{HLMV02,lO04b}, a local quantity defined in terms of first return times to a small neighbourhood of $x$ is considered, and the corresponding multifractal spectrum is studied.

Various asymptotic quantities associated to actions of Kleinian and Fuchsian groups are considered in~\cite{KS04,KS06,MV07}.

A notion of multifractal rigidity is studied in~\cite{BPS97b,lB08,BS08}, based on the idea that under certain circumstances, a finite number of multifractal spectra can serve as complete invariants within a particular class of dynamical systems.  That is, just as two Bernoulli shifts are isomorphic (as measure-preserving dynamical systems) if and only if they have the same measure-theoretic entropy, so it is expected that for broader classes of dynamical systems, there is a finite set of multifractal spectra $\FFF_1, \dots, \FFF_n$ with the property that two systems from the class are isomorphic (in the appropriate sense) if and only if their spectra $\FFF_1, \dots, \FFF_n$ are equal.

\section{Examples and applications}\label{sec:examples}

We describe various classes of systems and multifractal spectra to which the results in Section~\ref{sec:results} can be applied.  Many of these systems and spectra have already been studied in the papers surveyed in Section~\ref{sec:survey}.  In cases where previous work has used the thermodynamic approach, we obtain new results in the presence of phase transitions, which occurs for non-H\"older potentials and for non-uniformly hyperbolic maps.  In cases where previous work has used the saturation approach, we obtain new results in the form of the strengthened conditional variational principle~\eqref{eqn:cvp3}.

\subsection{Uniformly hyperbolic maps}\label{sec:unif-hyp}

\subsubsection{Entropy spectrum for Birkhoff averages}

Let $X$ be a compact metric space and $f\colon X\to X$ an expansive homeomorphism with specification.  (This includes the case when $f$ is an Axiom A diffeomorphism and $X$ is a basic set for $f$.)  Then it is shown in~\cite{rB75b} that every H\"older continuous potential $\ph$ has a unique equilibrium state.



Because $f$ is expansive, the entropy map is upper semi-continuous, and it is well-known that in this case, uniqueness of the equilibrium state is equivalent to differentiability of the pressure function in every direction.  That is, $\phi$ has a unique equilibrium state if and only if the map $q\mapsto P(\phi+q\psi)$ is differentiable at $q=0$ for every $\psi\in C(X)$.  (For example, see~\cite{dR78,gK98}.)  Thus $P$ is differentiable on the set of H\"older continuous potentials.

In particular, let $\ph\colon X\to \RR$ be H\"older continuous, and let $\psi=u\equiv 1$.  Then the multifractal spectrum associated to $(\ph,\psi,u)$ is the entropy spectrum for Birkhoff averages of $\ph$, and by Corollary~\ref{cor:birkhoff}, we have $I=\{\int \ph\,d\mu \mid \mu\in\Mf(X)\}$, $K_\alpha^\ph = \emptyset$ for all $\alpha\notin I$, and the spectrum is given by~\eqref{eqn:BEcvp} for all $\alpha\in \inter I$.  In particular, for every $\alpha\in \inter I$, there is a full measure $\nu_\alpha$, which is obtained as the equilibrium state for $q\ph$, where $q=q(\alpha)$ is such that $\frac{d}{dq}P(q\ph) = \alpha$.  This example is also given in~\cite{BS01}.

For general continuous potentials (not necessarily H\"older), it was shown in~\cite{TV03} that we have
\begin{equation}\label{eqn:TV}
\begin{aligned}
\htop(K_\alpha^\ph) &= \sup\left\{ h_\mu(f) \,\Big|\, \mu\in \Mf(X), \int \ph\,d\mu = \alpha\right\} \\
&= \inf \{ P(q\ph) - q\alpha \mid q\in \RR \}
\end{aligned}
\end{equation}
whenever $K_\alpha^\ph \neq \emptyset$.  (This was generalised in~\cite{PS07} to a broader class of systems satisfying a weaker version of the specification property.)  This establishes part, but not all, of~\eqref{eqn:BEcvp}.


To establish the final equality in~\eqref{eqn:BEcvp}, we need to know more about the thermodynamics of the potentials $q\ph$.  Thermodynamic properties of certain non-H\"older potentials were studied in~\cite{fH77,PZ06,hH08}.  For such potentials the function $T_0\colon q\mapsto P(q\ph)$ need not be differentiable; indeed, Pesin and Zhang show in~\cite{PZ06} that if $f$ is a uniformly piecewise expanding full-branched Markov map of the unit interval and $\ph\in C(X)$ is given by
\begin{equation}\label{eqn:nonHolder}
\ph(x) = \begin{cases}
-(1-\log x)^{-\beta} & x\in (0,1], \\ 0 & x=0 \end{cases}
\end{equation}
for some $0<\beta\leq 1$, then there exists $q_0>0$ at which $T_0$ is not differentiable.  Thus the results in~\cite{BS01} do not apply here, while Corollary~\ref{cor:birkhoff} shows that the formalism holds for all $\alpha\in \inter I \setminus (D^-T_0(q_0),D^+T_0(q_0))$; for all such $\alpha$ we have~\eqref{eqn:BEcvp}.

To extend the formalism to the interval $(D^-T_0(q_0),D^+T_0(q_0))$, on which the predicted spectrum is an affine function of $\alpha$, we need to use the last part of Corollary~\ref{cor:birkhoff}
---thus we need to find a family of compact $f$-invariant sets $X_n$ such that $P_{X_n}(q\ph) \to P_X(q\ph)$ for all $q$.

In Pesin and Zhang's example, the map $f$ is topologically (semi-)conjugate to the full shift $\Sigma_d^+$.  (It is conjugate on a total probability set, and hence by the variational principle, the pressure function for $(X,f)$ is entirely determined by the pressure function for $(\Sigma_d^+,\sigma)$.)  Writing $x_0$ for the sequence $000\dots$, we see that $\ph$ is H\"older on every compact subset of $\Sigma_d^+$ that does not contain $x_0$.  Let $X_n$ denote the set of sequences $x=x_0 x_1 \dots \in \Sigma_d^+$ in which the symbol $0$ never appears more than $n$ times consecutively.  Then $X_n$ is compact and $f$-invariant and $\ph$ is H\"older on $X_n$, so $q\mapsto P_{X_n}(q\ph)$ is differentiable for all $q$.  It remains only to show that $P_{X_n}(q\ph) \to P_X(q\ph)$.

\begin{proposition}\label{prop:enough-pressure}
Let $X_A$ be a transitive topological Markov chain on finitely many symbols, and $\{x_1,\dots,x_p\} \subset X_A$ a periodic orbit.  For $\eps>0$, let $X_\eps := \{x\in X_A \mid d(\sigma^n(x),\sigma^n(x_i)) \geq \eps \text{ for all } i,n\}$.  Then the $X_\eps$ are compact $\sigma$-invariant sets such that for every  $\phi\in C(X)$, we have $\lim_{\eps\to 0} P_{X_\eps}(\phi) = P_X(\phi)$.
\end{proposition}	


Proposition~\ref{prop:enough-pressure} lets us apply the last part of Corollary~\ref{cor:birkhoff} to obtain a new result for a class of non-H\"older potentials, which does not follow either from~\cite{BS01} or~\cite{TV03}.

\begin{theorem}[Entropy spectrum for Birkhoff averages of a non-H\"older potential]\label{thm:nonHolder}
Let $(X_A,\sigma)$ be a transitive topological Markov chain on finitely many symbols, and let $\ph\in C(X)$ be locally H\"older continuous everywhere except a single periodic orbit.  Then in addition to~\eqref{eqn:TV}, the level sets $K_\alpha^\ph = \{x \mid \frac 1n S_n \ph(x) \to \alpha \}$ satisfy
\[
\htop K_\alpha^\ph = \sup \{ h_\mu(f) \mid \mu \in \Mfe(K_\alpha^\ph) \}
\]
for every $\alpha \in \inter \{ \int \ph\,d\mu \mid \mu\in \Mf(X) \}$.
\end{theorem}

\subsubsection{Dimension spectra for conformal repellers}\label{sec:PW}

Let $M$ be a smooth manifold and $f\colon M\to M$ a $C^{1+\eps}$ map.  Let $X\subset M$ be compact and $f$-invariant ($f^{-1}X=X$).  Suppose that for every $x\in X$ there exists $a(x) > 1$ such that $Df(x)/a(x)\colon T_x M \to T_{f(x)}M$ is an isometry; then we say that $X$ is a conformal repeller for $f$.  Note that because $f$ is $C^{1+\eps}$, the function $a\colon X\to (1,\infty)$ is $\eps$-H\"older, and so is the function $\log a$.

As before, every H\"older continuous potential $\phi$ on $X$ has a unique equilibrium state, which is a Gibbs measure.  Because $f$ is expansive, the entropy map is upper semi-continuous, so every continuous potential has an equilibrium state and the pressure function $P$ is differentiable on the space of H\"older continuous potentials.  This thermodynamic information is all we need to apply the results from Section~\ref{sec:results}.

For example, the dimension spectrum for Lyapunov exponents can be computed via Corollary~\ref{cor:Lyapunov}, which shows that with $I=\{\lambda(\mu) \mid \mu\in \Mf(X)\}$, we have $K_\alpha^\lambda=\emptyset$ whenever $\alpha\notin I$, while the spectrum satisfies~\eqref{eqn:LDcvp} for all $\alpha\in \inter I$.  This was originally shown in~\cite{hW99}.  

For other potentials $\ph\in C(X)$ not equal to $\log a$, the dimension spectrum for Birkhoff averages of $\ph$ is a mixed spectrum of the sort studied in~\cite{BS01} under the assumption that $P$ is differentiable on $\spn\{\ph,\log a\}$.  Theorem~\ref{thm:full-measure} recovers one of the main results of that paper, namely that writing $I=\{\int \ph\,d\mu \mid \mu\in\Mf(X)\}$, when this differentiability occurs we have
\begin{equation}\label{eqn:BDcvp}
\begin{aligned}
\dim_H (K_\alpha^\ph) &= \inf \{ T_\alpha(q) \mid q\in \RR \} \\
&= \sup \left\{ \frac{h_\mu(f)}{\lambda(\mu)} \,\Big|\, \mu\in \Mf(X), \int \ph\,d\mu = \alpha \right\} \\
&= \sup \left\{ \frac{h_\mu(f)}{\lambda(\mu)} \,\Big|\, \mu\in \Mfe(K_\alpha^\ph) \right\}
\end{aligned}
\end{equation}
for every $\alpha\in \inter I$, where $T_\alpha(q)$ is defined by $P(q(\ph-\alpha) -T_\alpha(q) \log a) = 0$.

Even without full differentiability, Proposition~\ref{prop:enough-pressure} can sometimes be used to go beyond the results in~\cite{BS01}.  If $\ph$ is not H\"older continuous, then the pressure function is not necessarily differentiable on the span of $\{\ph,\log a\}$; thus the following result is new.

\begin{theorem}[Dimension spectrum for Birkhoff averages of a non-H\"older potential]
Let $f\colon M\to M$ be $C^{1+\eps}$ and let $X$ be a conformal repeller for $f$.  Let $\ph\in C(X)$ be H\"older continuous everywhere except for a single periodic orbit.  Then with $I$ as above, we have $K_\alpha^\ph = \emptyset$ for all $\alpha\notin I$, while $\dim_H (K_\alpha^\ph)$ is given by~\eqref{eqn:BDcvp} for every $\alpha\in \inter I$.
\end{theorem}
\begin{proof}
By Proposition~\ref{prop:enough-pressure}, there exist compact invariant sets $X_n\subset X$ on which $P_{X_n}\colon (q,t) \mapsto P_{X_n}(q\ph - t\log a)$ is differentiable, and furthermore $P_{X_n}(\phi) \to P_X (\phi)$ for all $\phi\in C(X)$, so Corollary~\ref{cor:phase} shows that~\eqref{eqn:BDcvp} holds for all $\alpha\in \RR$.
\end{proof}

Returning to the case where $\phi\in C(X)$ is H\"older continuous, we recall that its unique equilibrium state $\mu$ is a Gibbs measure, and so Corollary~\ref{cor:dimspec} recovers one of the main results of~\cite{PW97}, namely that for all $\alpha\in \inter I$, the dimension spectrum for pointwise dimensions of $\mu$ is given by~\eqref{eqn:DDcvp}.

\subsection{Maps with parabolic fixed points}\label{sec:MP}

Let $X$ be a compact metric space and $f\colon X\to X$ a conformal map; let $f'\in C(X)$ be the conformal derivative of $f$.  Suppose $f$ has a finite generating Markov partition whose boundary is given zero weight by every invariant measure $\mu$ with $\int \log f'\,d\mu > 0$.

Suppose that there is a fixed point $p=f(p)$ with $f'(p)=1$, that $f'(x)>1$ for all $x\neq P$, and that $f'$ is locally H\"older continuous on $X\setminus \{p\}$.

Finally, suppose that near $p$, the pair $(X,f)$ looks like a repeller for a conformal map on a manifold.  More precisely, suppose that there is a neighbourhood $U \subset X$ of $p$ and a Riemannian manifold $M$ into which $U$ embeds isometrically via a map $\pi$.  Let $f_*$ denote the pushforward of $f$ under $\pi$, and suppose that $f_*$ extends to a conformal map $\tilde f\colon Y\to M$ on some geodesically complete neighbourhood $Y$ of $\pi(p)$, with the property that for all large $n$ and small $\eps$, the function $\overline{B(p,n,\eps)} \to \RR$ given by $x\mapsto \|(\tilde f^n)'(x)\|$ attains its supremum on $\di B(p,n,\eps)$.

\begin{proposition}\label{prop:uhasP}
If $f$ satisfies the above conditions, then the function $u=\log f'$ satisfies \PP, and so $\dim_u$ is a well-defined Carath\'eodory dimension characteristic.
\end{proposition}

\begin{remark}
Proposition~\ref{prop:uhasP} gives the first example of a setting in which $u$-dimension is a Carath\'eodory dimension characteristic despite the fact that $u$ is not uniformly positive.
\end{remark}

\begin{example}
Define $f\colon [0,1]\to [0,1]$ by
\[
f(x) = \begin{cases}
\frac x{1-x} & \text{ for } 0\leq x \leq \frac 12, \\
\frac {2x-1}x & \text{ for } \frac 12 < x \leq 1.
\end{cases}
\]
Identifying the endpoints of the interval, we can think of $f$ as a smooth map of the unit circle with a single parabolic fixed point.  We have $f'(x) = (1-x)^{-2}$ for $x\in [0,1/2]$ and $f'(x)=x^{-2}$ for $x\in (1/2,1]$; thus $f$ satisfies the conditions above.
\end{example}

Using the same argument as in Proposition~\ref{prop:enough-pressure}, there exist compact $f$-invariant sets $X_n \subset X\setminus \{p\}$ such that $\lim_{n\to\infty} P_{X_n}(\phi) = P_X(\phi)$ for all $\phi\in C(X)$.  Furthermore, we may assume that $X_n$ has a finite Markov partition whose boundary is invisible to invariant measures (since $u$ is uniformly positive on $X_n$), and so for all H\"older continuous $\phi$, the pressure function on $X_n$ is differentiable on the span of $\{\phi,u\}$.  (Note that we may not have differentiablity of the pressure function on $X$ itself due to the parabolic fixed point, at which $u$ may not be H\"older.)

\subsubsection{Dimension spectrum for Lyapunov exponents}

It is shown in~\cite{GR09} that for $I = \{ \lambda(\mu) \mid \mu\in \Mf(X) \}$ and $\alpha\in \inter I$, we have
\begin{equation}\label{eqn:GR}
\begin{aligned}
\dim_H (K_\alpha^\lambda) &= \frac 1\alpha \inf \{ P(q\log f') - q\alpha \mid q\in \RR \} \\
&= \sup \left\{ \frac{h_\mu(f)}{\lambda(\mu)} \,\Big|\, \mu\in \Mf(X), \lambda(\mu) = \alpha \right\}.
\end{aligned}
\end{equation}
In light of the above discussion, Corollary~\ref{cor:Lyapunov} implies not only~\eqref{eqn:GR}, but also the following new result.

\begin{theorem}[Dimension spectrum for Lyapunov exponents of a non-uniformly expanding map]\label{thm:dimLyap}
Let $(X,f)$ be as above.  Then in addition to~\eqref{eqn:GR}, the level sets $K_\alpha^\lambda$ satisfy
\[
\dim_H (K_\alpha^\lambda) = \sup \left\{ \frac{h_\mu(f)}{\lambda(\mu)} \,\Big|\, \mu \in \Mfe(K_\alpha^\lambda)\right\}.
\]
\end{theorem}

\begin{remark}
Although Theorem~\ref{thm:dimLyap} is stronger than the corresponding results in~\cite{GR09}, we observe that our thermodynamic techniques do not suffice to prove the results in~\cite{GR09} regarding the irregular set or points with zero lower Lyapunov exponent.
\end{remark}

\subsubsection{Dimension spectrum for Birkhoff averages}

Moving to the dimension spectrum for Birkhoff averages, it is shown in~\cite{JJOP08} that for $\ph\in C(X)$, $I = \{\int \ph\,d\mu \mid \mu\in \Mf(X)\}$, and $\alpha\in \inter I$, we have
\begin{equation}\label{eqn:JJOP}
\begin{aligned}
\dim_H (K_\alpha^\ph) &= \inf \{ T_\alpha(q) \mid q\in \RR \} \\
&= \sup \left\{ \frac{h_\mu(f)}{\lambda(\mu)} \,\Big|\, \mu\in \Mf(X), \int \ph\,d\mu = \alpha \right\}.
\end{aligned}
\end{equation}
Under the stronger assumption that $\ph$ is H\"older continuous, we can use Theorem~\ref{thm:universal} and Corollary~\ref{cor:phase} to say this and more.  Those results apply to the level sets $\hat K_\alpha$, and so we first observe that $\hat K_\alpha = K_\alpha$ for nearly every value of $\alpha$.

\begin{lemma}\label{lem:llzero}
Let $f$ be a conformal expanding map with a parabolic fixed point $p$ satisfying the conditions above, and let $\ph\in C(X)$.  Then for every $\alpha \neq \ph(p)$, we have $K_\alpha \subset \hat X$, whence $\hat K_\alpha = K_\alpha$.
\end{lemma}

This implies the following new result.

\begin{theorem}[Dimension spectrum for Birkhoff averages of a non-uniformly expanding map]
Let $(X,f)$ be as above and let $\ph\in C(X)$ be H\"older continuous.   Then in addition to~\eqref{eqn:JJOP}, the level sets $K_\alpha^\ph = \{x \mid \frac 1n S_n \ph(x) \to \alpha \}$ satisfy
\[
\dim_H (K_\alpha^\ph) = \sup \left\{ \frac{h_\mu(f)}{\lambda(\mu)} \,\Big|\, \mu\in \Mfe(K_\alpha^\ph) \right\}
\]
for every $\alpha \in \inter \{ \int \ph\,d\mu \mid \mu\in \Mf(X) \}$.
\end{theorem}

\begin{remark}
As in the previous section, we observe that our thermodynamic techniques do not yield the stronger results in~\cite{JJOP08} on the irregular set, and more work is required to describe the precise shape of the spectrum and the value of $\dim_H K_{\ph(p)}$ as is done there.
\end{remark}

\subsubsection{Dimension spectrum for pointwise dimensions}

Finally, we consider the dimension spectrum for pointwise dimensions.  Let $\mu$ be a weak Gibbs measure for some potential $\phi\in C(X)$.  (For the existence of such measures, see~\cite{mK01,JR09}.)  Writing $\ph = P(\phi) - \phi$, it is shown in~\cite{JR09} that for $I=\{\frac{\int \ph\,d\mu}{\lambda(\mu)} \mid \mu \in \Mf(X), \lambda(\mu)>0\}$ and $\alpha\in \inter I$, we have
\begin{equation}\label{eqn:JR}
\begin{aligned}
\dim_H (K_\alpha^{d_\mu}) &= \inf \{ T_0(q) - q\alpha \mid q\in \RR \} \\
&= \sup \left\{ \frac{h_\nu(f)}{\lambda(\nu)} \,\Big|\, \nu\in \Mf(X), \frac{\int \ph\,d\nu}{\lambda(\nu)} = \alpha \right\}.
\end{aligned}
\end{equation}

If we assume in addition that $\phi$ is H\"older continuous and $P(\phi) > \phi(p)$, then we can obtain this and more from Corollary~\ref{cor:dimspec}.

\begin{lemma}\label{lem:llzero2}
With $X$, $f$, $\phi$, and $\mu$ as above, we have~\eqref{eqn:uphi} for all $x\in X$.
\end{lemma}

Using Lemma~\ref{lem:llzero2}, we can apply Corollary~\ref{cor:dimspec} to obtain the following new result.

\begin{theorem}[Dimension spectrum for pointwise dimensions of a non-uniformly expanding map]
Let $X$, $f$, $\phi$, and $\mu$ be as above, and suppose that $\dim_u X < \infty$.  Then in addition to~\eqref{eqn:JR}, the level sets $K_\alpha^{d_\mu} = \{x \mid d_\mu(x)=\alpha\}$ satisfy
\[
\dim_H (K_\alpha^{d_\mu})
= \sup \left\{ \frac{h_\nu(f)}{\lambda(\nu)} \,\Big|\, \nu \in \Mfe(K_\alpha^{d_\mu})\right\}
\]
for every $\alpha\in \inter I$.
\end{theorem}
\begin{proof}
Any measure $\nu\in \Mfe(X)$ with $\int u\,d\nu = 0$ must have $\nu(\{p\})=1$, whence $h_\nu(f)=0$.  This verifies the hypotheses of Corollary~\ref{cor:dimspec}.
\end{proof}


\subsection{Other non-uniformly hyperbolic maps}

\subsubsection{Maps with critical points---continuous potentials}

For interval maps with critical points (unimodal and multimodal maps), existence and uniqueness of equilibrium states for a certain class of potentials were established in~\cite{BT08}.  In particular, let $\HHH$ denote the collection of topologically mixing $C^2$ interval maps $f\colon [0,1]\to[0,1]$ with hyperbolically repelling periodic points and non-flat critical points and fix $f\in \HHH$.

It was shown by Blokh that any continuous topologically mixing interval map has the specification property~\cite{aB83,jB97}.  Thus for any continuous potential $\ph\in C([0,1])$, the saturation results in~\cite{PS07} can be applied (note that the results in~\cite{TV03} do not apply because $f$ is not expansive).  For the level sets $K_\alpha^\ph = \{x \mid \frac 1n S_n\ph(x) \to \alpha\}$, this yields
\begin{equation}\label{eqn:PfSu}
\begin{aligned}
\htop(K_\alpha^\ph) &= \sup\left\{ h_\mu(f) \,\Big|\, \mu\in \Mf(X), \int \ph\,d\mu = \alpha\right\} \\
&= \inf \{ P(q\ph) - q\alpha \mid q\in \RR \}
\end{aligned}
\end{equation}
for all $\alpha \in \inter \{ \int \ph\,d\mu \mid \mu\in \Mf(X) \}$.  Once again, we can strengthen~\eqref{eqn:PfSu} given appropriate thermodynamic information.

Let $\ph\colon [0,1]\to \RR$ be a H\"older continuous potential such that
\begin{equation}\label{eqn:BR}
\sup \ph - \inf \ph < \htop (f).
\end{equation}
It is shown in~\cite{BT08} that there exists a unique equilibrium state for $\ph$, and so by upper semi-continuity of the entropy map~\cite[Lemma 2.3]{BK98}, there exists $q_0>1$ such that $T_0\colon q\mapsto P(q\ph)$ is differentiable on $(-q_0,q_0)$ and we can apply Corollary~\ref{cor:birkhoff} to show that the entropy spectrum for Birkhoff averages of $\ph$ is given by~\eqref{eqn:BEcvp} on $(D^+T_0(-q_0),D^-T_0(q_0))$.  Thus we have the following new result.

\begin{theorem}[Entropy spectrum for Birkhoff averages of a map with critical points]\label{thm:unimodal}
Fix $f\in \HHH$ and let $\ph\in C([0,1])$ be H\"older continuous such that~\eqref{eqn:BR} holds.  Then for every $\alpha \in (D^+T_0(-q_0),D^-T_0(q_0))$, in addition to~\eqref{eqn:PfSu}, the level set $K_\alpha^\ph$ satisfies
\[
\htop(K_\alpha^\ph) = \sup \{ h_\mu(f) \mid \mu\in \Mfe(K_\alpha^\ph) \}.
\]
\end{theorem}

\subsubsection{Maps with critical points---discontinuous potentials}

For maps with critical points, the geometric potential $\ph(x) = -\log |f'(x)|$ is no longer continuous or even bounded, and so the results presented here do not apply directly to the multifractal analysis of the spectra associated to this potential or to the various dimension spectra associated to $f$.  Nevertheless, the thermodynamics of the potentials $t\ph$ have recently been studied in~\cite{PS08,BT09,IT09b}, and it is reasonable to expect that these thermodynamic results can be used to obtain multifractal information.  This has been done in~\cite{mT09,IT09} using inducing schemes; a more abstract deduction of (weaker) multifractal results using methods along the lines of this paper has been given in the preprint~\cite{vC10c}.

\subsubsection{Viana's examples}

The existence and uniqueness of equilibrium states for a broad class of non-uniformly expanding maps in higher dimensions was studied by Oliveira and Viana~\cite{OV08} and by Varandas and Viana~\cite{VV08}.  The multifractal properties of these systems do not appear to have been studied, nor does the question of whether or not these systems (which may have contracting regions) satisfy specification or any of its variants.

We briefly describe the systems studied in~\cite{VV08}.  Let $M$ be a compact connected manifold and let $f\colon M\to M$ be a local homeomorphism.  Writing $d$ for distance on $M$, suppose there is a bounded function $L(x)$ such that every $x\in M$ has a neighbourhood $U_x\ni x$ on which $f_x = f|_{U_x} \colon U_x \to f(U_x)$ is invertible, with $d(f(y),f(z)) \geq L(x)^{-1} d(y,z)$ for all $y,z\in U_x$.  Thus if $L(x) <1$, then $f$ is expanding at $x$, while if $L(x)\geq 1$, then $L$ controls how much contraction can happen near $x$.

Assuming every point has finitely many preimages, we write $\deg_x(f) = \# f^{-1}(x)$.  If level sets for the degree are closed, then is it shown in~\cite{VV08} that up to considering some iterate $f^N$ of $f$, we can assume that $\deg_x(f)\geq e^{\htop(f)}$ for all $x$.

Finally, suppose there are $\sigma>1$ and an open region $\AAA\subset M$ such that
\begin{enumerate}[(H1)]
\item  $L(x)\leq \sigma^{-1}$ for all $x\in M\setminus \AAA$, and $\sup_{x\in \AAA} L(x)$ is close to $1$ (see~\cite{VV08} for precise conditions);
\item  there exists a covering $\PPP$ of $M$ by domains of injectivity for $f$ such that $\AAA$ can be covered by $r<e^{\htop(f)}$ elements of $\PPP$.
\end{enumerate}

Thus $f$ is uniformly expanding outside of $\AAA$, and does not display too much contraction inside $\AAA$; furthermore, since there are at least $e^{\htop(f)}$ preimages of any given point $x$, and only $r$ of these can lie in covering of $\AAA$ by elements of $\PPP$, every point has at least one preimage in the expanding region.  (See~\cite{VV08} for examples of specific systems satisfying these conditions.)

Now suppose $\ph\in C(M)$ is H\"older continuous and satisfies $\sup\ph - \inf \ph < \htop(f) - \log r$.  Then it is proved in~\cite{VV08} that there exists a unique equilibrium state for $\ph$.  By upper semi-continuity of the entropy map, this implies differentiability of the pressure function at $\ph$.

As in the discussion before Theorem~\ref{thm:unimodal}, there exists $q_0>1$ such that $T_0\colon q\mapsto P(q\ph)$ is differentiable on $(-q_0,q_0)$ and we can apply Corollary~\ref{cor:birkhoff} to give the following new result.

\begin{theorem}
Let $M$, $f$, and $\ph$ be as above.  Then we have
\begin{align*}
\htop (K_\alpha^\ph) &= \lF(\alpha) = \uF(\alpha) \\
&= \inf \{ T_0(q) - q\alpha \mid q\in \RR \} \\
&= \sup \left\{ h_\mu(f) \,\Big|\, \mu\in \Mf(X), \int \ph\,d\mu = \alpha \right\} \\
&= \sup \{ h_\mu(f) \mid \mu \in \Mfe(K_\alpha^\ph)\}
\end{align*}
for every $\alpha \in (D^+T_0(-q_0),D^-T_0(q_0))$.
\end{theorem}

\section{Proofs}\label{sec:pfs}

\subsection{Proof of Proposition~\ref{prop:ugeq0}}

It suffices to show that $\int u\,d\mu \geq 0$ whenever $\mu\in \Mf(X)$ is ergodic.  Suppose $\mu\in \Mfe(X)$ is such that $\int u\,d\mu < 0$, and let $x$ be generic for $\mu$, so $\frac 1n S_n u(x) \to \int u\,d\mu < 0$.  Let $\eps>0$ be such that $S_n u(x) < -n\eps$ for all sufficiently large $n$, and let $\delta>0$ be such that $|u(z) - u(y)| < \eps$ whenever $d(y,z) < \delta$.  Then if $y\in X$ and $n\in \NN$ are such that $n$ is large and $x\in B(y,n,\delta)$, we have $u(f^k(y)) \leq u(f^k(x)) + \eps$ for all $0\leq k< n$, whence $S_n u(y) \leq S_n u(x) + n\eps \leq 0$.  This contradicts property \PP.

\subsection{Proof of Proposition~\ref{prop:bowen}}

We will use the following lemma a number of times.

\begin{lemma}\label{lem:bound-pressure}
Given $f\colon X\to X$, $\eta, \phi\in C(X)$, and $Z\subset X$, suppose there exist $\alpha, \beta \in \RR$ such that
\[
\alpha \leq \llim_{n\to\infty} \frac 1n S_n\phi(x) \leq \ulim_{n\to\infty} \frac 1n S_n\phi(x) \leq \beta
\]
for every $x\in Z$.  Then
\begin{equation}\label{eqn:cone}
P_Z(\eta) + \alpha t \leq P_Z(\eta + t\phi) \leq P_Z(\eta) + \beta t
\end{equation}
for all $t>0$.
\end{lemma}
\begin{proof}This follows from the proof of Proposition 5.3 in \cite{vC10a}.
\end{proof}

Now fix $Z\subset X$ as in the statement of Proposition~\ref{prop:bowen} and observe that
\[
m_P(Z,0,-tu,N,\delta) = \inf_{\PPP(Z,N,\delta)} \sum_{(x_i,n_i)} e^{S_{n_i}(-tu)} =
m_u(Z,t,N,\delta)
\]
for all $t\in \RR$, $N\in \NN$, and $\delta>0$.  Thus $m_P(Z,0,-tu,\delta) = m_u(Z,t,\delta)$.

Let $t^* = \inf\{t\in \RR \mid P_Z(-tu) \leq 0\}$.  Given $t < t^*$, we have $P_Z(-tu,\delta) > 0$ for all sufficiently small $\delta>0$, and hence
\[
m_u(Z,t,\delta) = m_P(Z,0,-tu,\delta) = +\infty,
\]
which implies that $\dim_u(Z,\delta) \geq t$.  Since $t<t^*$ was arbitrary, this shows that $\dim_u Z \geq t^*$.

Now consider the sets $Z_m = \{ x\in Z \mid \llim_{n\to\infty} \frac 1n S_n u(x) \geq \frac 1m \}$.  By the assumption that $\llim \frac 1n S_n u(x) > 0$ for all $x\in Z$, we have $\bigcup_{m\geq 1} Z_m = Z$, and hence by~\eqref{eqn:dimustable}, to show that $\dim_u Z \leq t^*$ it suffices to show that $\dim_u Z_m \leq t^*$ for all $m$.

To this end, fix $m\in \NN$ and $t>t^*$.  By Lemma~\ref{lem:bound-pressure}, we have
\[
P_{Z_m}(-tu) \leq P_{Z_m}(-t^* u) - \frac{t-t^*}m \leq P_Z(-t^* u) - \frac{t-t^*}m < 0,
\]
and hence for all sufficiently small $\delta$ we have
\[
m_u(Z_n,t,\delta) = m_P(Z,0,-tu,\delta) = 0,
\]
whence $\dim_u Z_m \leq t$.  Since $t>t^*$ was arbitrary, this implies that $\dim_u Z_m \leq t^*$, and we are done.

\begin{remark}
Observe that if $Z=\emptyset$, then $P_Z(-tu)=-\infty$ for all $t\in \RR$, which agrees with the convention that $\dim_u Z = -\infty$.  If $Z\neq \emptyset$, then there exists $x\in Z$ with $\llim \frac 1n S_n u(x) = \alpha > 0$, whence for all $t<0$ we have $P_Z(-tu) \geq P_{\{x\}}(-tu) \geq |t| \alpha > 0$, and so it suffices to consider $t\geq 0$ in Bowen's equation.
\end{remark}

\subsection{Proof of Proposition~\ref{prop:dimumu}}

Let $G_\mu$ be the set of generic points for $\mu$, so $\frac 1n S_n u(x) \to \int u\,d\mu$ for every $x\in G_\mu$, and $\htop(G_\mu) = h_\mu(f)$ by~\cite[Theorem 3]{rB73}.  Using Lemma~\ref{lem:bound-pressure}, we obtain
\[
P_{G_\mu}(-tu) = P_{G_\mu}(0) - t\int u\,d\mu = \htop (G_\mu) - t\int u\,d\mu = h_\mu(f) - t\int u\,d\mu.
\]
This is equal to $0$ if and only if $t=h_\mu(f)/\int u\,d\mu$, so it follows from the Birkhoff ergodic theorem and Proposition~\ref{prop:bowen} that $\dim_u \mu \leq \dim_u G_\mu = h_\mu(f) / \int u\,d\mu$.

For the other inequality, we observe that if $Z\subset X$ is such that $\mu(Z) = 1$, then writing $Z' = \{x\in Z\mid \frac 1n S_n u(x) \to \int u\,d\mu\}$, we have $\mu(Z')=1$, whence $\htop(Z') \geq h_\mu(f)$, and so using Lemma~\ref{lem:bound-pressure} again, we have
\[
P_{Z'}(-tu) = \htop(Z') - t\int u\,d\mu \geq h_\mu(f) - t\int u\,d\mu,
\]
whence by Proposition~\ref{prop:bowen} we have $\dim_u Z \geq \dim_u Z' \geq h_\mu(f)/\int u\,d\mu$.

\subsection{Proof of Proposition~\ref{prop:Sgeq0}}

\begin{lemma}\label{lem:no-coming-back}
If $\phi\in C(X)$ is such that $P(\phi)\leq 0$, then $P(\lambda\phi) \leq 0$ for all $\lambda\geq 1$.
\end{lemma}
\begin{proof}
For every $\mu\in \Mf(X)$, the fact that $P(\phi) \leq 0$ implies that $h_\mu(f) + \int \phi\,d\mu \leq 0$.  In particular, $\int \phi\,d\mu \leq 0$, hence for all $\lambda\geq 1$,
\[
h_\mu(f) + \int \lambda\phi\,d\mu = \left(h_\mu(f) + \int \phi\,d\mu\right) + (\lambda-1) \int \phi\,d\mu \leq 0 + 0 = 0.
\]
It follows that $P(\lambda\phi) = \sup_{\mu\in \Mf(X)} (h_\mu(f) + \int \lambda\phi\,d\mu) \leq 0$.
\end{proof}

Now suppose $\SSS(\alpha) < 0$.  Then there exist $q,t\in \RR$ with $t<0$ such that $P(q(\ph - \alpha \psi) - tu) \leq 0$.  By Lemma~\ref{lem:no-coming-back}, we have $R_\alpha(\lambda q,\lambda t) \leq 0$ for every $\lambda \geq 1$, and so $T_\alpha(\lambda q) \leq \lambda t$.  In particular, since $t < 0$, this implies that $\SSS(\alpha) = \inf T_\alpha = -\infty$, and so $\alpha\notin I$.

\subsection{Proof of Theorem~\ref{thm:universal}}

Fix $\alpha \in \RR$ and $t < \hat\FFF(\alpha)$.  (If $\hat\FFF(\alpha)=-\infty$, then the desired inequality is automatic.)  We claim that
\begin{equation}\label{eqn:Raqt}
R_\alpha(q,t) = P(q(\ph - \alpha\psi) - tu) > 0,
\end{equation}
for every $q$, and hence $T_\alpha(q) \geq t$.  This will in turn imply that
\[
\SSS(\alpha) = \inf \{ T_\alpha(q) \mid q\in \RR \} \geq t,
\]
and since $t < \hat\FFF(\alpha)$ was arbitrary, this suffices.  Thus it remains only to prove~\eqref{eqn:Raqt}.

\begin{lemma}\label{lem:ratiotoother}
If $\frac{S_n\ph(x)}{S_n\psi(x)} \to \alpha$, then $\frac 1nS_n(\ph-\alpha\psi)\to 0$.
\end{lemma}
\begin{proof}
$|\frac 1nS_n(\ph-\alpha\psi)| = |\frac{S_n\ph(x)}{S_n\psi(x)} - \alpha|\cdot |\frac 1n S_n\psi(x)| \leq|\frac{S_n\ph(x)}{S_n\psi(x)} - \alpha| \sup \psi$.
\end{proof}

Using Bowen's equation (Proposition~\ref{prop:bowen}), the fact that $t < \hat\FFF(\alpha) = \dim_u \hat K_\alpha$ and $\llim_{n\to\infty} \frac 1n S_n u(x) > 0$ for all $x\in \hat K_\alpha$ implies that $P_{\hat K_\alpha}(-tu) > 0$.  Furthermore, for every $q\in\RR$, Lemmas~\ref{lem:bound-pressure} and~\ref{lem:ratiotoother} imply that $P_{\hat K_\alpha} (q(\ph-\alpha\psi) - tu) = P_{\hat K_\alpha}(-tu)$, which establishes~\eqref{eqn:Raqt}.

\subsection{Proof of Theorem~\ref{thm:full-measure}}

The key tool is Ruelle's formula for the derivative of pressure.

\begin{proposition}\label{prop:ruelle}
Fix $\eta,\phi\in C(X)$.  If the function $q\mapsto P(\eta + q\phi)$ is differentiable at $q=0$, and if in addition $\nu$ is an equilibrium state for $\eta$, then
\begin{equation}\label{eqn:derivative}
\frac{d}{dq} P(\eta + q\phi)|_{q=0} = \int \phi \,d \nu.
\end{equation}
\end{proposition}
\begin{proof}
Write $g(q) = P(\eta + q\phi)$.  Then for all $q\in\RR$, we have
\begin{multline*}
g(q) = \sup_\mu \left\{ h(\mu) + \int \eta \,d\mu + \int q\phi \,d\mu \right\} \\
\geq h(\nu) + \int \eta \,d\nu + q \int \phi \,d\nu 
= P(\eta) + q \int \phi \,d\nu,
\end{multline*}
whence $g(q) - g(0) = g(q) - P(\eta) \geq q \int \phi \,d\nu$.  In particular, for $q>0$, we get $\frac 1q(g(q) - g(0)) \geq \int \phi \,d\nu$, and hence $g'(q)\geq \int \phi \,d\nu$ (recall that differentiability of $g$ was one of the hypotheses), while for $q<0$, we have $\frac 1q(g(q)-g(0)) \leq \int \phi \,d\nu$, and hence $g'(q)\leq \int \phi \,d\nu$, which establishes equality.
\end{proof}

In the setting of Theorem~\ref{thm:full-measure}, we observe that $R_\alpha(q,T_\alpha(q)) = 0$ for all $q$, and so
\begin{align*}
0 &= \frac d{dq} P(q(\ph - \alpha\psi) - T_\alpha(q) u) \\
&= \frac d{dq} P(q(\ph - \alpha \psi) - \SSS(\alpha)u) + \frac d{dq} T_\alpha(q) \frac d{dt} P(q(\ph - \alpha\psi) - tu)|_{t=\SSS(\alpha)}.
\end{align*}
The fact that $q$ realises the infimum in~\eqref{eqn:S} implies that $\frac d{dq}T_\alpha(q) = 0$, and hence by Proposition~\ref{prop:ruelle}, we have $\int (\ph - \alpha\psi) \,d\nu_\alpha = 0$.  This yields $\alpha = \int \ph\,d\nu_\alpha / \int \psi\,d\nu_\alpha$ (recall we assume that $\int \psi\,d\nu_\alpha > 0$), and ergodicity of $\nu_\alpha$ implies that $\nu_\alpha (K_\alpha) = 1$, whence $\FFF(\alpha) \geq \dim_u \nu_\alpha$.

Because $\nu_\alpha$ is an equilibrium state for $q(\ph - \alpha\psi) -\SSS(\alpha)u$, we have
\[
0 = P(q(\ph - \alpha\psi) - \SSS(\alpha) u) = h(\nu_\alpha) - \SSS(\alpha) \int u\,d\nu_\alpha,
\]
and so (using the assumption that $\int u\,d\nu_\alpha > 0$),
\[
\SSS(\alpha) = \frac{h(\nu_\alpha)}{\int u\,d\nu_\alpha} = \dim_u \nu_\alpha \leq \FFF(\alpha).
\]

\subsection{Proof of Proposition~\ref{prop:Ihatint}}

First we show that $\hat I \subset \inter I$ without any hypotheses regarding $\int\psi\,d\mu$ or $h_\mu(f)$.  To this end, fix $\alpha\in \hat I$, so $0\in \inter \KKK(\alpha)$.  Thus there exist $\eta>0$ and $\mu^\pm \in \Mf(X)$ such that $\int(\ph-\alpha\psi)\,d\mu^\pm =\pm\eta$.  Choosing $\eps>0$ sufficiently small, we have $\int(\ph-\alpha'\psi)\,d\mu^+ > 0$ and $\int(\ph-\alpha'\psi)\,d\mu^- <0$ for all $\alpha'\in (\alpha-\eps,\alpha+\eps)$.  It follows that $0\in \inter \KKK(\alpha')$ for each such $\alpha'$, so $\hat I$ is open.

Furthermore, given $\alpha\in \hat I$ and $\eta>0$, $\mu^\pm\in \Mf(X)$ as before, we have
\[
R_\alpha(q,0) = P(q(\ph-\alpha\psi)) \geq \max \left(h_{\mu^\pm}(f) + \int q(\ph-\alpha\psi)\,d\mu^\pm\right)
\geq |q| \eta > 0
\]
for all $q\neq 0$, whence $T_\alpha(q) > 0$.  It follows that $\SSS(\alpha) \geq 0$, and so $\alpha\in I$.  Now since $\hat I$ is open and contained in $I$, we have $\hat I \subset \inter I$.

It remains to show that $\inter I \subset \hat I$.  To this end, let $\gamma>0$ be such that $\int \psi\,d\mu \geq \gamma h_\mu(f)$ for all $\mu\in \Mfe(X)$, and fix $\alpha\notin \hat I$.  Then $0\notin \inter \KKK(\alpha)$, and since $\KKK(\alpha)$ is an interval, we either have $\KKK(\alpha)\subset [0,\infty)$ or $\KKK(\alpha)\subset (-\infty,0]$.

First consider the case $\KKK(\alpha) \subset [0,\infty)$.  Then $\int (\ph-\alpha\psi)\,d\mu\geq 0$ for all $\mu\in \Mfe(X)$, and so for every $\eps>0$ and $q>0$, we have
\begin{align*}
R_{\alpha+\eps}(q,t) &= P(q(\ph-(\alpha+\eps)\psi) - tu) \\
&= \sup_{\mu\in \Mfe(X)} \left( h_\mu(f) + q\int (\ph-\alpha\psi)\,d\mu - q\eps\int \psi\,d\mu - t\int u\,d\mu \right) \\
&\leq \sup_{\mu\in \Mfe(X)} (h_\mu(f)(1 - q\eps\gamma) + |t|\sup u ).
\end{align*}
For sufficiently large $q$ and for $t<0$ with $|t|$ sufficiently small, this gives $R_{\alpha+\eps}(q,t) \leq 0$, whence $\SSS(\alpha+\eps) \leq T_{\alpha+\eps}(q) \leq t < 0$, and so $\alpha+\eps\notin I$.  Since $\eps>0$ was arbitrary, we have $\alpha\notin \inter I$.  The proof in the case $\KKK(\alpha)\subset (-\infty,0]$ is similar, and we see that $\inter I = \hat I$.

\subsection{Proof of Theorem~\ref{thm:Scvp}}

Let $\Mfa(X)$ be as in~\eqref{eqn:Mfa}, and write $\tilde\SSS(\alpha) := \sup \{\frac {h_\mu(f)}{\int u\,d\mu} \mid \mu \in \Mfa(X)\}$.  We begin by showing that $\SSS(\alpha) \geq \tilde\SSS(\alpha)$ for every $\alpha$: to this end, fix $t < \tilde\SSS(\alpha)$.  Then there exists $\mu\in \Mfa(X)$ such that $h_\mu(f) > t\int u\,d\mu$, and therefore for every $q\in \RR$ we have
\[
P(q(\ph - \alpha\psi) - tu) \geq h(\mu) + q\int(\ph-\alpha\psi)\,d\mu - t\int u\,d\mu > 0,
\]
so $T_\alpha(q) > t$.  This implies that $\SSS(\alpha) \geq t$, and since $t < \tilde\SSS(\alpha)$ was arbitrary, one direction is done.

It remains only to show that $\tilde\SSS(\alpha) \geq \SSS(\alpha)$ when $\alpha\in \hat I$.  Fix $\alpha\in \hat I$, and for each $q,t\in \RR$, write $\phi_{q,t} := q(\ph-\alpha\psi) - tu$.  Consider the collections of measures
\begin{align*}
\MMM^+ &:= \left\{\mu\in \Mfa(X) \,\Big|\, \int(\ph-\alpha\psi)\,d\mu \geq 0 \right\}, \\
\MMM^- &:= \left\{\mu\in \Mfa(X) \,\Big|\, \int(\ph-\alpha\psi)\,d\mu \leq 0 \right\},
\end{align*}
and define two functions $g^\pm\colon \RR^2 \to \RR$ by
\[
g^\pm(q,t) := \sup\left\{h_\mu(f) + \int \phi_{q,t}\,d\mu \,\Big|\, \mu\in \MMM^\pm \right\}.
\]
Observe that $P(\phi_{q,t}) = \max(g^+(q,t),g^-(q,t))$, and so writing
\begin{equation}\label{eqn:Epm}
\begin{aligned}
E^\pm &:= \{(q,t) \in \RR^2 \mid P(\phi_{q,t}) = g^\pm(q,t)\}\\
&= \{(q,t) \in \RR^2 \mid g^\pm(q,t) \geq g^\mp(q,t)\},
\end{aligned}
\end{equation}
we have $E^+ \cup E^- = \RR^2$.  (The sets $E^+$ and $E^-$ need not be disjoint.)

It follows from the definitions that $g^-(q,t)$ is non-increasing in $q$ and $g^+(q,t)$ is non-decreasing in $q$.  Using the second characterisation in~\eqref{eqn:Epm}, this implies that given any $(q,t)\in E^-$, we have $(q',t)\in E^-$ and $P(\phi_{q',t}) \geq P(\phi_{q,t})$ for all $q'<q$.  Similarly, given any $(q,t)\in E^+$, we have $(q',t)\in E^+$ and $P(\phi_{q',t}) \leq P(\phi_{q,t})$ for all $q'>q$.

Now fix $t < \SSS(\alpha)$ and let $E^\pm(t) := \{q\in \RR \mid (q,t)\in E^\pm\}$.  Because $0\in \inter \KKK(\alpha)$ (which follows from $\alpha\in \hat I)$, there exists $\mu \in \Mf(X)$ with  $\int(\ph-\alpha\psi)\,d\mu > 0$, whence $P(\phi_{q,t}) > P(\phi_{0,t})$ for sufficiently large values of $q$; for such $q$ we have $q\in E^+(t)$, and so $E^+(t)\neq \emptyset$.  A similar argument shows that $E^-(t) \neq \emptyset$.

Because both $E^-(t)$ and $E^+(t)$ are non-empty, there exists $q_0\in \overline{E^-(t)} \cap \overline{E^+(t)}$.  If $q_0\in E^-(t)$, then we have $q_0+\gamma\in E^+(t)$ for all $\gamma>0$, while if $q_0\in E^+(t)$, then $q_0-\gamma\in E^-(t)$ for all $\gamma>0$.  Without loss of generality, we assume the former case; the proof in the latter case is similar.

Now we have $q^-:=q_0\in E^-(t)$ such that $q^+\in E^+(t)$ for all $q^+>q^-$.  Let $\delta := \inf \{ P(\phi_{q,t})\mid |q-q^-|\leq 1\}$, and observe that $\delta>0$.  Choose $q^+ > q^-$ so that $q^+ - q^- < \frac \delta2 \|\ph-\alpha\psi\|^{-1}$.  Then there exist $\mu^\pm\in \MMM^\pm$ such that
\[
h_{\mu^\pm}(f) + \int \phi_{q^\pm}\,d\mu^\pm \geq \frac \delta2.
\]
Choose $s\in [0,1]$ such that $\mu:= s\mu^- + (1-s)\mu^+$ satisfies $\int (\ph-\alpha\psi)\,d\mu=0$.  Then we have
\begin{align*}
h_\mu(f) - &t\int u\,d\mu = h_\mu(f) + \int \phi_{q^-,t}\,d\mu \\
&= sh_{\mu^-} + (1-s)h_{\mu^+} + s\int \phi_{q^-,t}\,d\mu^- + (1-s)\int\phi_{q^-,t}\,d\mu^+ \\
&\geq \frac\delta2 + (1-s)\int (\phi_{q^-,t} - \phi_{q^+,t})\,d\mu^+ \\
&\geq \frac\delta2 - (q^+-q^-) \|\ph - \alpha\psi\| > 0,
\end{align*}
where the final inequality follows from the choice of $q^+$.  Let $q=q(\alpha)$ be such that $T_\alpha(q) = \SSS(\alpha)$---such a $q$ exists because $0\in \inter \KKK(\alpha)$---and observe that 
\begin{multline*}
h_\mu(f) - \SSS(\alpha) \int u\,d\mu = h_\mu(f) + \int \phi_{q,\SSS(\alpha)}\,d\mu \\
\leq P(\phi_{q,\SSS(\alpha)}) = 0 < h_\mu(f) - t\int u\,d\mu,
\end{multline*}
whence $\int u\,d\mu > 0$, so $\mu\in \Mfa(X)$.  Finally, the inequality $h_\mu(f) - t\int u\,d\mu > 0$ implies that $t < \frac{h_\mu(f)}{\int u\,d\mu}$, and since $t < \SSS(\alpha)$ was arbitrary, this yields $\tilde\SSS(\alpha) \geq \SSS(\alpha)$.

\subsection{Proof of Proposition~\ref{prop:Qfinite}}

The proposition is a direct consequence of the following lemma, which applies even when $u\neq \psi$.

\begin{lemma}\label{lem:toinfinity}
For every $\lambda\geq 1$ and $\alpha,q\in \RR$, we have $T_\alpha(\lambda q) \leq \lambda T_\alpha(q)$.  In particular, if $q\in \RR$ is such that $T_\alpha(q) < \infty$, then $T_\alpha(\lambda q)<\infty$ for all $\lambda \geq 1$.
\end{lemma}
\begin{proof}
$T_\alpha(q) < \infty$ if and only if there exists $t < \infty$ such that $P(q(\ph-\alpha\psi) - tu) \leq 0$.  It follows from Lemma~\ref{lem:no-coming-back} that in this case $P(\lambda q(\ph-\alpha\psi) - \lambda t u) \leq 0$ for all $\lambda \geq 1$, whence $T_\alpha(\lambda q) \leq \lambda T_\alpha(q) < \infty$.
\end{proof}

\subsection{Proof of Theorem~\ref{thm:concave-hull}}

To show that $\SSS$ is the concave hull of $\hat\FFF$ on $\hat I$, it suffices to show that
\[
T_0(q) = \sup \{ \hat\FFF(\alpha) + q\alpha \mid \alpha\in \RR \}
\]
for every $q\in \inter Q$, since then $\SSS$ is the double Legendre transform of $\hat\FFF$ on $\hat I$, which is the concave hull.  Thus we prove the following two inequalities for every $q\in \inter Q$:
\begin{align}
\label{eqn:Tgeq}
T_0(q) &\geq \hat\FFF(\alpha) + q\alpha \text{ for all } \alpha\in \RR, \\
\label{eqn:Tleq}
T_0(q) &\leq \sup \{ \hat\FFF(\alpha) + q\alpha \mid \alpha \in \RR \}.
\end{align}
It follows from Theorem~\ref{thm:universal} that
\[
\hat\FFF(\alpha)\leq\SSS(\alpha)\leq T_0(q) - q\alpha
\]
for every $q,\alpha$, which shows~\eqref{eqn:Tgeq}.

Turning our attention to~\eqref{eqn:Tleq}, we see that if $\nu \in \Mfe(X)$ is such that $h_\nu(f) + q\int \ph\,d\nu - t\int u\,d\nu >0$ and $\int u\,d\nu=0$, then for all $t'>t$, we have
\begin{align*}
P(q\ph - t'u) &\geq h_\nu(f) + q\int \ph\,d\nu - t'\int u\,d\nu  \\
&= h_\nu(f) + q\int \ph\,d\nu - t\int u \,d\nu > 0
\end{align*}
whence $T_0(q) = +\infty$.  Thus given $q\in \inter Q$, $t < T_0(q)$, and potentials of the form $q\ph - tu$, the supremum in the variational principle may be taken over measures with $\int u\,d\nu>0$.

For such a $q$ and $t$, let $\nu\in \Mfe(X)$ be such that $h(\nu) + q\int\ph\,d\nu - t\int u\,d\nu > 0$ and $\int u\,d\nu>0$, and let $\alpha$ be such that $\nu(K_\alpha)=1$.  Then we have
\[
(\dim_u \nu + q\alpha - t) u(\nu) > 0,
\]
which immediately yields $t < \FFF(\alpha) - q\alpha$.  Since $t < T_0(q)$ was arbitrary, this proves~\eqref{eqn:Tleq}.

\subsection{Proof of Proposition~\ref{prop:Ihatint2}}

Let $\mu\in \Mfe(X)$.  By the hypothesis, we have $h_\mu(f) = 0$ whenever $\int u\,d\mu=0$, and so it remains only to consider measures with $\int u\,d\mu > 0$.  For such a measure~\eqref{eqn:dimumu} and Proposition~\ref{prop:dimumu} give $\dim_u X \geq \dim_u \mu = h_\mu(f) / \int u\,d\mu$, which suffices.

\subsection{Proof of Theorem~\ref{thm:universal2}}

The proof of Theorem~\ref{thm:universal2} comes in two parts.  First we compare the fine and coarse multifractal spectra, showing that $\FFF \leq \lF$, the first inequality in~\eqref{eqn:FFS}; then we compare the coarse spectrum with the thermodynamically predicted spectrum, showing that $\uF \leq \SSS$, the third inequality in~\eqref{eqn:FFS}.  (The second inequality in~\eqref{eqn:FFS} is immediate.)

\subsubsection{Comparison of fine and coarse spectra}

We begin by expressing the level sets $K_\alpha$ from~\eqref{eqn:level2} in terms of the approximate level sets $G_n(U)$ from~\eqref{eqn:approx}.  We see that
\[
K_\alpha = \bigcap_{U\ni \alpha} \bigcup_{N\in\NN} \bigcap_{n\geq N} G_n(U).
\]
Thus for each neighbourhood $U\ni \alpha$, we consider the sets
\begin{equation}\label{eqn:Faen}
\begin{aligned}
\tG_N(U) &= \bigcap_{n\geq N} G_n(U) = \left\{ x\in X\,\Big|\, \frac 1n S_n\ph(x) \in U \text{ for all } n\geq N\right\} \\
\tG(U) &= \bigcup_{N\in\NN} \tG_N(U),
\end{aligned}
\end{equation}
for which we have $K_\alpha = \bigcap_{U\ni\alpha} \tG(U)$.  In particular, we see that for every neighbourhood $U\ni \alpha$,
\begin{equation}\label{eqn:dimuKG}
\htop K_\alpha \leq \htop \tG(U) = \sup_N \htop \tG_N(U) \leq \sup_N \lhtop \tG_N(U),
\end{equation}
where $\lhtop$ is the lower capacity topological entropy.  (See Section~\ref{sec:dim} for the notion of a capacity.)  Mimicking the definition of $\Lambda_n^\delta(U)$, let
\[
\tilde\Lambda_n^\delta(U,N) := \inf \{ \#E \mid B(E,n,\delta) \supset \tG_N(U) \},
\]
and observe that since $G_n(U) \supset \tG_N(U)$ for all $n\geq N$, we have
\[
\Lambda_n^\delta(U) \geq \tilde\Lambda_n^\delta(U,N)
\]
for all $n\geq N$.  In particular,
\[
\lim_{\delta\to 0} \llim_{n\to\infty} \frac 1n \log \Lambda_n^\delta(U) \geq
\lim_{\delta\to 0} \llim_{n\to\infty} \frac 1n \log \tilde\Lambda_n^\delta(U,N) 
= \lhtop \tG_N(U)
\]
for every $N$, and so~\eqref{eqn:dimuKG} yields
\[
\lim_{\delta\to 0} \llim_{n\to\infty} \frac 1n \log \Lambda_n^\delta(U) \geq \htop F(U) \geq \htop K_\alpha = \FFF(\alpha).
\]
This holds for all neighbourhoods $U\ni \alpha$, and so $\lF(\alpha) \geq \FFF(\alpha)$.

\subsubsection{Comparison of coarse and predicted spectra}

Now we prove the final inequality in~\eqref{eqn:FFS}.  Fix $\alpha\in \RR$ and let $U$ be any neighbourhood of $\alpha$.  Then for every $t < \uF(\alpha)$ and $\eps>0$ there exists $\delta>0$ and a sequence $n_k\to \infty$ such that $\Lambda_{n_k}^\delta(U) \geq e^{n_k t}$ for all $k$ and $|\ph(x)-\ph(y)| < \eps$ whenever $d(x,y)<\delta$.

Let $E\subset X$ be an $(n_k,\delta)$-spanning set for $X$ and let
\[
E' = \{x\in E \mid B(x,n,\delta) \cap G_n(U) \neq  \emptyset \}.
\]
Then from the definition of $\Lambda_{n_k}^\delta(U)$, we have $\# E' \geq e^{n_k t}$, and so
\[
\sum_{x\in E} e^{S_n(q\ph)(x)} \geq
\sum_{x\in E'} e^{S_n(q\ph)(x)} \geq 
e^{n_k t} e^{n_k q (\inf U - \eps)} =
 e^{n_k (t + q(\inf U-\eps))},
\]
which gives $P(q\ph) \geq t + q(\inf U-\eps)$.  Since $t < \uF(\alpha)$, $U\ni \alpha$, and $\eps>0$ were arbitrary, this implies that $P(q\ph) \geq \uF(\alpha) + q\alpha$.

\subsection{Proof of Corollary~\ref{cor:birkhoff}}

We begin by showing~\eqref{eqn:Ibirk}.  The first equality is straightforward; recall that $\alpha\in I$ if and only if $T_\alpha(q)\geq 0$ for all $q$, which in this case is equivalent to the inequality $P(q\ph) \geq q\alpha$ for all $q$.  It follows from properties of convex functions that this is true if and only if
\[
\alpha\in \left[\lim_{q\to-\infty} D^-T_0(q), \lim_{q\to\infty} D^+T_0(q)\right]
= \overline{\bigcup_{q\in\RR} J(q)}.
\]
For the second equality in~\eqref{eqn:Ibirk}, we first observe that if $\mu\in \Mf(X)$ is such that $\int \ph\,d\mu=\alpha$, then we have $P(q\ph) \geq h_\mu(f) + \int q\ph\,d\mu \geq q\alpha$ for all $q\in \RR$.  The other inclusion goes as follows.  Given $q\in \RR$, let $\alpha^+ = D^+T_0(q)$.  For $\gamma>0$, let $\mu_\gamma$ be an equilibrium state for $(q+\gamma)\ph$, and let $\mu\in \Mf(X)$ be a weak* limit of $\mu_\gamma$ as $\gamma\to 0$.  By upper semi-continuity, $\mu$ is an equilibrium state for $q\ph$, and moreover, $\int \ph\,d\mu = \lim_{\gamma\to 0} \int \ph\,d\mu_\gamma = \alpha^+$ since $\int \ph\,d\mu_\gamma \in J(q+\gamma)$ for all $\gamma$.  Because $\alpha^+ \in \di J(q)$, there exists an ergodic component $\nu^+$ of $\mu$ which is also an equilibrium state for $q\ph$ and which has $\int \ph\,d\nu^+ = \alpha^+$.

A similar argument produces an ergodic measure $\nu^-$ which is an equilibrium state for $q\ph$ and which has $\int \ph\,d\nu^- = \alpha^-$.  Taking convex combinations of $\nu^\pm$ shows that $\{\int \ph\,d\mu \mid \mu\in \Mf(X)\} \supset J(q)$ for every $q\in \RR$, which completes the characterisation of $I$.  (Note that the set of values taken by the integral is closed since it is a continuous image of a compact set.)

The completeness result $K_\alpha = \emptyset$ for all $\alpha\notin I$ follows from Theorem~\ref{thm:universal}.  When $\alpha$ is such that $J(q) = \{\alpha\}$ for some $q\in \RR$, the first two lines of~\eqref{eqn:BEcvp} follow from Theorems~\ref{thm:full-measure} and~\ref{thm:universal2}.  When $\alpha\in \di J(q)$ for some $q$, they follow from Theorem~\ref{thm:universal2} and the argument above that produces an ergodic equilibrium state $\nu$ for $q\ph$ with $\int \ph\,d\nu = \alpha$.

The third line of~\eqref{eqn:BEcvp} follows from Theorem~\ref{thm:Scvp}, and the fourth follows since we proved the first two lines by constructing a full measure.  Finally, the last statement of Corollary~\ref{cor:birkhoff} follows from Corollary~\ref{cor:phase}.

\subsection{Proof of Corollary~\ref{cor:Lyapunov}}

Once we establish the form of $\SSS$ and $I$, the proof here is virtually identical to the proof of Corollary~\ref{cor:birkhoff}.  Given $\alpha>0$ and $t,q\in \RR$, we see that the following are equivalent:
\begin{enumerate}
\item $t < T_\alpha(q)$;
\item $P(qu-q\alpha -tu) > 0$;
\item $P((q-t)u) - q\alpha  > 0$;
\item $P(q'u) - q'\alpha > t\alpha$, where $q'=q-t$;
\item $t < \frac 1\alpha (T_0(q') - q'\alpha)$.
\end{enumerate}
Thus since $t < \SSS(\alpha)$ if and only if $t<T_\alpha(q)$ for every $q$, which is equivalent to $t<\frac 1\alpha (T_0(q') - q'\alpha)$ for every $q'$, we see that $\SSS$ has the form claimed.  It follows immediately that for $\alpha>0$, we have $\SSS(\alpha)\geq 0$ if and only if $\inf \{T_0(q) - q\alpha \mid q\in \RR\}\geq 0$, and then~\eqref{eqn:ILyap} follows from~\eqref{eqn:Ibirk}, with the possible exception of $\alpha=0$.

For $\alpha=0$, we have $T_\alpha(q) = T_0(q)$ and so $\SSS(\alpha) = \inf_q T_\alpha(q) = \inf_q T_0(q)$.  Thus $0\in I$ if and only if $T_0(q) \geq 0$ for all $q\in \RR$.  This is the same criterion as in Corollary~\ref{cor:birkhoff}, which completes the proof of~\eqref{eqn:ILyap}.

Applying the results from~\cite{vC10a} as in the discussion, the rest of Corollary~\ref{cor:Lyapunov} follows as Corollary~\ref{cor:birkhoff} did.

\subsection{Proof of Corollary~\ref{cor:dimspec}}

Let $\tilde X = \{x\in X \mid \text{\eqref{eqn:uphi} holds} \}$, and observe that under the assumptions of Corollary~\ref{cor:dimspec}, we have $\tilde X=X$ and $X^\star = \hat X$.  In fact, we prove the corollary under the slightly weaker assumptions that writing $K_\alpha = K_\alpha^{d_\mu} \cup K_\alpha^{\ph,u}$, we have $K_\alpha \subset \tilde X$ and $K_\alpha\cap \hat X \subset X^\star$ for every $0\leq \alpha<\infty$.

\begin{proposition}\label{prop:dimspec}
Let $X$ be a compact metric space and $f\colon X\to X$ be conformal (without critical points or singularities).  Suppose $\mu\in \MMM(X)$ and $\phi\in C(X)$ satisfy~\eqref{eqn:wkgibbs}, and let $\psi = u = \log a$, $\ph = \phi-P(\phi)$.  Suppose $u$ satisfies \PP.  Then for every $0\leq \alpha < \infty$, we have
\begin{equation}\label{eqn:leveltilde}
\tilde X \cap K_\alpha^{\ph,u} = \tilde X \cap K_\alpha^{d_\mu} \subset \hat X,
\end{equation}
and in particular,
\begin{equation}\label{eqn:levelstar}
\dim_u (K_\alpha^{\ph,u} \cap X^\star) = \dim_H (K_\alpha^{d_\mu} \cap X^\star).
\end{equation} 
\end{proposition}
\begin{proof}
Fix $x\in \tilde X$ and observe that by~\cite[Lemma 6.1]{vC10a}, for every $\eps>0$ there exist $\delta_0=\delta_0(\eps)>0$ and $\eta=\eta(x,\eps)>0$ such that for every $n\in \NN$ and $0<\delta<\delta_0$,
\begin{equation}\label{eqn:bdddist}
B\left(x,\eta\delta e^{-(S_n u(x) + n\eps)} \right)
\subset B(x,n,\delta)
\subset B\left(x,\delta e^{-(S_n u(x) - n\eps)} \right).
\end{equation}
Furthermore, by the weak Gibbs property~\eqref{eqn:wkgibbs}, for every $\gamma>0$ we may assume by taking $\delta>0$ sufficiently small that
\begin{equation}\label{eqn:wkgibbs2}
|S_n \ph(x) + \log \mu(B(x,n,\delta))| \leq n\gamma
\end{equation}
for all large $n$.  Define $r_1,r_2\colon \NN\to (0,\delta)$ by
\[
r_1(n) = \delta e^{-(S_n u(x) - n\eps)}, \qquad
r_2(n) = \eta \delta e^{-(S_n u(x) + n\eps)},
\]
so that~\eqref{eqn:bdddist} becomes
\begin{equation}\label{eqn:bdddist2}
B(x,r_2(n)) \subset B(x,n,\delta) \subset B(x,r_1(n)).
\end{equation}
Observe that because $u$ is bounded away from $\pm\infty$, there exists $\xi>0$ such that $\frac{r_i(n+1)}{r_i(n)} \in [\xi^{-1},\xi]$ for all $n$.  In particular, $\ld_\mu(x)$ and $\ud_\mu(x)$ can be computed by considering balls with radius $r_i(n)$.

It follows from~\eqref{eqn:bdddist2} and~\eqref{eqn:wkgibbs2} that
\[
\log \mu(B(x,r_2(n))) \leq \log \mu(B(x,n,\delta)) \leq -(S_n\ph(x) - n\gamma),
\]
whence
\begin{equation}\label{eqn:lower-estimate}
\frac{\log \mu(B(x,r_2(n)))}{\log r_2(n)} \geq \frac{-(S_n\ph(x) - n\gamma)}{-(S_n u(x) + n\eps) + \log (\delta \eta)}.
\end{equation}
If $\llim \frac 1n S_n u(x) = 0$, then there exists $n_k\to \infty$ such that $\frac 1{n_k} S_{n_k} u(x) \to 0$.  Since $x\in \tilde X$, this implies that $\llim \frac 1{n_k} S_{n_k} \ph(x) > 0$; in particular, we can choose $\gamma$ small enough that $S_{n_k} \ph(x) \geq 2n_k \gamma$ for all $k$, whence~\eqref{eqn:lower-estimate} gives
\[
\frac{\log \mu(B(x,r_2(n_k)))}{\log r_2(n_k)} \geq \frac {n_k \gamma}{S_{n_k} u(x) + n_k\eps - \log (\delta\eta)}.
\]
Taking the limit as $k\to\infty$ gives $\ud_\mu(x) \geq \gamma/\eps$, and since $\eps>0$ was arbitrary, this shows that $d_\mu(x)\neq \alpha$ for all $0\leq \alpha < \infty$, whence $K_\alpha^{d_\mu} \cap \tilde X \subset \hat X$.

Similarly, given $x\in \tilde X \setminus \hat X$, we have $n_k$ such that $\frac{S_{n_k} \ph(x)}{S_{n_k} u(x)} \to \infty$, whence $\frac{S_n\ph(x)}{S_nu(x)} \not\to \alpha$ for any $0\leq \alpha<\infty$; thus $K_\alpha^{\ph,u} \cap \tilde X \subset \hat X$.

Now to establish~\eqref{eqn:leveltilde}, we can restrict our attention to the case where $x\in \hat X$---that is, $\llim \frac 1n S_n u(x) > 0$.

Given $x\in \hat X$, fix $\xi>0$ and suppose $\eps,\gamma$ are small enough such that $\xi (\frac 1n S_n u(x)) > \max\{\eps,\gamma\}$ for all large $n$.  Then~\eqref{eqn:lower-estimate} gives
\[
\frac{\log \mu(B(x,r_2(n)))}{\log r_2(n)} \geq \frac{S_n\ph(x) - \xi(S_n u(x)))}{(S_n u(x))(1+\xi) - \log (\delta \eta)},
\]
and since $\xi$ can be arbitrarily small, we have $\ld_\mu(x) \geq \llim \frac{S_n \ph(x)}{S_n u(x)}$.

Similarly, it follows from~\ref{eqn:bdddist2} and~\ref{eqn:wkgibbs2} that
\[
\log \mu(B(x,r_1(n))) \geq \log \mu(B(x,n,\delta)) \geq -(S_n\ph(x) + n\gamma),
\]
whence with $\xi$ as before, we have
\[
\frac{\log \mu(B(x,r_1(n)))}{\log r_1(n)} \leq \frac{-(S_n\ph(x) + n\gamma)}{-(S_n u(x) - n\eps) + \log \delta}
\leq \frac{S_n\ph(x) + \xi(S_n u(x)))}{(S_n u(x))(1-\xi) - \log \delta },
\]
and so $\ud_\mu(x) \leq \ulim \frac{S_n \ph(x)}{S_n u(x)}$.  This finishes the proof of~\eqref{eqn:leveltilde}.  As pointed out in Section~\ref{sec:Lyapunov}, we have $\dim_u Z = \dim_H Z$ for all $Z\subset X^\star$, and~\eqref{eqn:levelstar} follows.
\end{proof}

With Proposition~\ref{prop:dimspec} in hand, the proof of Corollary~\ref{cor:dimspec} follows exactly the same lines as the proof of Corollary~\ref{cor:birkhoff}, using the results from Sections~\ref{sec:mainresult}--\ref{sec:concave}.  The only difference is that we need the hypotheses that $\dim_u X<\infty$ and $\int u\,d\mu>0$ for every $\mu\in\Mf(X)$ with $h_\mu(f)>0$ in order to apply Proposition~\ref{prop:Ihatint2} and show that Theorem~\ref{thm:Scvp} applies on all of $\inter I$.  For the form of $I$, it suffices to observe that the following statements are equivalent:
\begin{enumerate}
\item $\alpha \in \overline{\bigcup_{q\in \RR} J(q)}$;
\item $T_0(q) \geq q\alpha$ for all $q\in \RR$;
\item $P(q\ph - t\log f') > 0$ for every $q\in \RR$ and $t<q\alpha$;
\item $P(q\ph - (q\alpha + t')\log f') > 0$ for every $q\in \RR$ and $t'<0$;
\item $\SSS(\alpha) \geq 0$.
\end{enumerate}

\subsection{Proofs of results in examples}\label{sec:egpfs}

\begin{proof}[Proof of Proposition~\ref{prop:enough-pressure}]
By passing to an iterate of $\sigma$ and recoding, we can assume that $X_A \setminus \bigcup_n X_n = \{ 0^\infty \}$ is a single fixed point.  Once again, passing to a higher iterate and recoding allows us to assume that the transitions $0\to 1$ and $1\to 0$ are allowed by the transition matrix $A$.

Let $E(k) = \{0,1,\dots,d-1\}^k$, and for every $w\in E(k)$, let $[w] = \{x\in X_A \mid x_i=w_i \text{ for all } 1\leq i\leq d \}$.  Let $E_A(k) = \{w\in E(k) \mid [w]\neq \emptyset\}$, and let $E_n(k)$ be the set of words in $E_A(k)$ which do not contain more than $n$ consecutive appearances of the symbol $0$.  Observe that it suffices to prove the proposition in the case $X_n = \{x\in X_A \mid [x_1\dots x_k] \in E_n(k) \text{ for all } k \}$.

Define $\phi_k\colon E_A(k)\to \RR$ by $\phi_k(w) = \sup \{ S_n \phi(x) \mid x\in [w] \}$.  Consider the partition sums
\[
\Lambda_k(X_A,\phi) = \sum_{w \in E_A(k)} e^{\phi_k(w)}, \qquad\qquad
\Lambda_k(X_n,\phi) = \sum_{w\in E_n(k)} e^{\phi_k(w)},
\]
and observe that $P_{X_A}(\phi) = \lim_{k\to\infty} \frac 1k \log\Lambda_k(X_A,\phi)$, and similarly for $X_n$.  Now define $\pi\colon E_A(k) \to E_n(k)$ by replacing every $n$th $0$ with a $1$ in each string of more than $n$ consecutive $0$s.  We estimate the cardinalities $|\pi^{-1}(w)|$.

Fix $w\in E_n(k)$ and let $J=\{j\in \{1,\dots,k\} \mid w_i = 0 \text{ for all } j-n \leq i < j \}$.  Then we have $|J|\leq k/n$, and furthermore, if $\pi(w') = w$, then $w_i = w_i'$ for all $i\notin J$, while for each $j\in J$ we have $w_i' \{0,1\}$.  Thus $|\pi^{-1}(w)| \leq 2^{k/n}$.

Furthermore, if $\pi(w')=w$, then $w'$ and $w$ disagree in at most $k/n$ indices, which are separated by at least $n$.  Let $m_1,\dots,m_p \geq n$ be such that
\[
\{m_1 + m_2 + \cdots + m_j \mid 1\leq j \leq p\} = \{i\in \{1,2,\dots,k\} \mid w_i' \neq w_i \},
\]
and observe that $p \leq k/n$.  Writing $V_\ell(\phi) = \sup \{\phi(x)-\phi(y) \mid x_i=y_i \text{ for all } 1\leq i\leq \ell \}$, we have
\[
|\phi_k(w') - \phi_k(w)| \leq \sum_{i=1}^{p} \sum_{\ell=1}^{m_i} V_\ell(\phi) \leq \frac kn \sum_{\ell=1}^n V_\ell(\phi),
\]
where the last inequality follows from the fact that $V_\ell(\phi) \leq V_{\ell'}(\phi)$ for all $\ell'\geq \ell$.  Now we can compare the partition sums, and hence the pressures:
\begin{align*}
\Lambda_k(X_A,\phi) &= \sum_{w\in E_n(k)} \sum_{w'\in \pi^{-1}(w)} e^{\phi_k(w')} \\
&\leq 2^{k/n} e^{(k/n)\sum_{\ell=1}^n V_\ell(\phi)} \Lambda_k(X_n,\phi), \\
\frac 1k\log\Lambda_k(X_A,\phi) &\leq \frac 1n\left( \log 2 + \sum_{\ell=1}^n V_\ell(\phi) \right) + \frac 1k\log\Lambda_k(X_n,\phi).
\end{align*}
Take the limit as $k\to\infty$ and observe that since $\phi$ is uniformly continuous, we have $\lim_{\ell\to\infty} V_\ell(\phi)=0$, which implies that $\lim_{n\to\infty} \frac 1n \sum_{\ell=1}^n V_\ell(\phi) =0$.
\end{proof}

\begin{proof}[Proof of Proposition~\ref{prop:uhasP}]
First we show that for every $\eps>0$ there exist points $x_n\in B(p,n,\eps)$ such that $S_n u(x_n) \to \infty$.  Fix $0<\eps'<\eps$ and let $x_n\in \di B(p,n,\eps')$ be such that $S_n u(x_n) = \sup \{ S_n u(y) \mid y\in \overline{B(p,n,\eps')}\}$.

Let $\Delta \subset X\times X$ denote the diagonal $\{(x,x) \mid x\in X\}$, and observe that the function $X\times X\setminus \Delta \to \RR^+$ given by $(x,y) \mapsto d(f^n(x),f^n(y))/d(x,y)$ extends to a continuous function on $X\times X$ by $(x,x) \mapsto e^{S_n u(x)}$.  Because $X\times X$ is compact, this function is uniformly continuous, and so for every $\gamma>0$ there exists $\eta=\eta(\gamma,n)>0$ such that
\begin{equation}\label{eqn:almostright}
\left| \frac{d(f^n(x),f^n(y))}{d(x,y)} - e^{S_n u(x)} \right| < \gamma
\end{equation}
whenever $d(x,y) < \eta$.  Working now on the manifold $M$ and using geodesic completeness, let $y_0,\dots,y_k$ be such that $y_0=\pi(p)$, $y_k = \pi(x_n)$, and $d(y_i,y_{i+1}) < \eta$ for all $i$, and furthermore such that $\sum_i d(y_i,y_{i+1}) = d(x,p)$.  Choose $\gamma$ and $\eta$ such that~\eqref{eqn:almostright} extends to $\tilde f$ as well; then we have
\begin{align*}
d(f^n(x),p) &\leq \sum_{i=0}^{k-1} d(\tilde{f}^n(y_i),\tilde{f}^n(y_{i+1})) 
\leq \sum_{i=0}^{k-1} \left( e^{S_n u(y_i)} + \gamma\right) d(y_i,y_{i+1}) \\
&\leq \left(e^{S_n u(x_n)} + \gamma\right) d(x,p).
\end{align*}
Using the fact that $x\in \di B(p,n,\eps')$ and hence $d(f^n(x_n),p)=\eps'$, we see that $e^{S_n u(x_n)} \geq \eps' d(x_n,p)^{-1} - \gamma$.  Because the Markov partition is generating, we have $\diam B(p,n,\eps')\to 0$ as $n\to\infty$, whence $d(x_n,p)\to 0$ and hence $S_n u(x_n) \to \infty$.

It remains to show that \PP\ holds for each $\delta>0$.  Using the fact that $f$ has a generating Markov partition, pass to a refinement $\PPP$ with the property that $\diam \PPP < \delta$.  Let $\pi\colon \Sigma_A \to X$ be the coding map for the Markov partition $\PPP$, and suppose without loss of generality that $\pi(000\cdots)=p$.  

Let $L_N$ denote the collection of $N$-cylinders in $\Sigma_A$, and given an $N$-cylinder $[w]$, let $\phi_N(w) := \sup_{x\in [w]} S_N u(x)$.  Observe that to establish \PP, it suffices to show that $\lim_{N\to\infty} \inf_{[w]\in L_N} \phi_N(w)=\infty$, since then we may take as our collection $E_N$ the set $\{(x_w, N) \mid [w]\in L_N\}$, where $x_w$ is such that $\phi_N(w) = S_N u(x_w)$.

Let $\gamma := \inf \{u(x) \mid x_0 \neq 0 \}$, and observe that $\gamma > 0$.  Let $C_N := \phi_N(0^N)$, and observe that by the arguments in the first half of the proof, we have $C_N\to \infty$.  Now given $[w]\in L_N$, one of two things happens; either the word $w$ contains a string of at least $\sqrt{N}$ consecutive $0$s, in which case $\phi_N(w) \geq C_{\sqrt{N}}$, or $w$ contains at least $\sqrt{N}$ symbols that are not $0$, in which case $\phi_N(w) \geq \sqrt{N} \gamma$.  It follows that $\inf_{[w]\in L_N} \phi_N(w) \geq \min(\sqrt{N}\gamma, C_{\sqrt{N}})$, which suffices.
\end{proof}

\begin{proof}[Proof of Lemma~\ref{lem:llzero}]
Suppose $x\in X\setminus \hat X$, so $\llim_{n\to\infty} \frac 1n S_n u(x) = 0$.  Given $\delta>0$, let $\eps(\delta)>0$ be such that $u(x) \geq \eps(\delta)$ for all $x$ with $B(x,p) \geq \delta$.

By the hypothesis on $x$, there exists $n_k$ such that $\frac 1{n_k} S_{n_k} u(x) \to 0$.  Let $\delta_k\to 0$ be such that writing $\eps_k = \eps(\delta_k)$, we have $S_{n_k} u(x) \leq n_k \eps_k \delta_k$.  Let $q_n(x,\delta) = \#\{k\in \{0,1,\dots,n-1\} \mid d(f^k(x),p) \geq \delta$, and observe that $S_n u(x) \geq q_n(x,\delta) \eps(\delta)$.  In particular, we have
\[
q_{n_k}(x,\delta_k) \eps_k \leq S_{n_k} u(x) \leq n_k \eps_k \delta_k,
\]
whence $q_{n_k}(x,\delta_k) \leq \delta_k n_k$.  It follows immediately that $\lim_{k\to\infty} \frac 1{n_k} S_{n_k} \ph(x) = \ph(p)$, and so either $x\in K_{\ph(p)}$ or $\frac 1n S_n \ph(x)$ does not converge to a limit, in which case $x\notin K_\alpha$ for any $\alpha$.
\end{proof}

\begin{proof}[Proof of Lemma~\ref{lem:llzero2}]
Fix $x\in X$ and suppose there is a subsequence $n_k$ is such that $\frac 1{n_k} S_{n_k} u(x) \to 0$.  As in the proof of Lemma~\ref{lem:llzero}, we have $\frac 1{n_k} S_{n_k} \phi (x) \to \phi(p) < P(\phi)$, and so writing $\ph = P(\phi)-\phi$, we see that $\frac 1{n_k} S_{n_k} \ph(x) \to P(\phi)-\phi(p) < 0$.  It follows that $x$ satisfies~\eqref{eqn:uphi}.
\end{proof}

\section{Background material on Carath\'eodory dimension characteristics}\label{sec:dim}

We recall the general theory of Carath\'eodory dimension characteristics, following~\cite[Chapter 4]{yP98}, and then verify that the $u$-dimension defined by Barreira and Schmeling fits into this theory even if $u$ is not uniformly positive, so long as condition \PP\ is satisfied.  

Given a set $X$, a C-structure $\CCC = \{\SSS, \FFF, \xi, \eta, \psi \}$ on $X$ is specified by an arbitrary index set $\SSS$, a collection of subsets of $X$ indexed by $\SSS$ and denoted $\FFF = \{U_s \mid s\in \SSS\}$, and functions $\eta, \psi, \xi \colon \SSS\to \RR^+$ such that
\begin{description}
\item[A1]  $\emptyset \in \FFF$; if $U_s=\emptyset$ then $\eta(s) = \psi(s)=0$; if $U_s\neq \emptyset$ then $\eta(s)>0$ and $\psi(s)>0$;
\item[A2] for any $\delta>0$ there exists $\eps>0$ such that $\eta(s)\leq\delta$ for any $s\in \SSS$ with $\psi(s)\leq \eps$;
\item[A3] for any $\eps>0$ there exists a finite or countable subcollection $\GGG\subset \SSS$ that covers $X$ (that is, $\bigcup_{s\in \GGG} U_s \supset X$) and has $\psi(\GGG) := \sup\{\psi(s) \mid s\in \GGG\} \leq \eps$.
\end{description}
Note that the definition places no requirements on the function $\xi$.  Given a set $Z\subset X$ and $\eps>0$, consider the following collection of $\eps$-covers of $Z$:
\[
\PPP(Z,\eps) := \left\{\GGG\subset \SSS \,\Big|\, Z \subset \bigcup_{s\in \GGG} U_s \text{ and }
\psi(s) \leq \eps \text{ for all } s\in \GGG \right\}.
\]
Then condition A3 amounts to requiring that $\PPP(X,\eps)$ be non-empty for every $\eps>0$.  In fact, it is useful to require a slightly stronger condition: to wit, consider the collection
\[
\tilde\PPP(Z,\eps) := \left\{\GGG\subset \SSS \,\Big|\, Z \subset \bigcup_{s\in \GGG} U_s \text{ and }
\psi(s) = \eps \text{ for all } s\in \GGG \right\},
\]
and replace A3 by the condition that
\begin{description}
\item[A3$'$] there exist arbitrarily small $\eps>0$ such that $\tilde\PPP(X,\eps)\neq\emptyset$.
\end{description}
To obtain a Carath\'eodory dimension characteristic from a C-structure $\CCC$, one defines for every $Z\subset X$, $\alpha\in \RR$, and $\eps>0$ the numbers
\begin{align*}
M_\CCC(Z,\alpha,\eps) &= \inf \left\{ \sum_{s\in \GGG} \xi(s) \eta(s)^\alpha \,\Big|\, \GGG \in \PPP(Z,\eps) \right\}, \\
m_\CCC(Z,\alpha) &= \lim_{\eps\to 0} M_C(Z,\alpha,\eps),
\end{align*}
and considers the critical point of the function $\alpha\mapsto m_\CCC(Z,\alpha)$:
\begin{align*}
\dim_\CCC(Z) &= \inf \{ \alpha\in \RR \mid m_\CCC(Z,\alpha) = 0 \} \\
&= \sup \{\alpha\in \RR \mid m_\CCC(Z,\alpha) = \infty \}.
\end{align*}
Provided A3$'$ holds, one may also define the \emph{upper and lower $\CCC$-capacities} by
\begin{align*}
R_\CCC(Z,\alpha,\eps) &= \inf \left\{ \sum_{s\in \GGG} \xi(s) \eta(s)^\alpha \,\Big|\, \GGG \in \tilde\PPP(Z,\eps) \right\}, \\
\ur_\CCC(Z,\alpha) &= \ulim_{\eps\to 0} R_\CCC(Z,\alpha,\eps), \\
\lr_\CCC(Z,\alpha) &= \llim_{\eps\to 0} R_\CCC(Z,\alpha,\eps), \\
\uCap_\CCC(Z) &= \inf\{\alpha\in \RR \mid \ur_\CCC(Z,\alpha) = 0 \}, \\
\lCap_\CCC(Z) &= \inf\{\alpha\in \RR \mid \lr_\CCC(Z,\alpha) = 0 \},
\end{align*}
where the limit in the definition of $\ur_\CCC$ and $\lr_\CCC$ is taken over those values of $\eps$ for which $\tilde\PPP(Z,\eps)\neq\emptyset$.

\begin{remark}
In~\cite{yP98} condition A3$'$ asks for $\tilde\PPP(Z,\eps)$ to be non-empty for \emph{every} sufficiently small $\eps>0$.  However, there are many interesting C-structures (including a number considered in~\cite{yP98}) for which $\psi$ takes discrete values, and hence one only has the version of A3$'$ stated here.
\end{remark}

If in addition to the above conditions one has
\begin{description}
\item[A4] $\eta(s_1) = \eta(s_2)$ for any $s_1,s_2\in \SSS$ with $\psi(s_1) = \psi(s_2)$,
\end{description}
then defining $\eta(\eps)$ by $\eta(\eps) := \eta(s)$ for any $s$ with $\psi(s)=\eps$, and setting
\[
\Lambda_\CCC(Z,\eps) = \inf \left\{ \sum_{s\in \GGG} \xi(s) \,\Big|\, \GGG \in \tilde\PPP(Z,\eps) \right\},
\]
it can be shown that
\[
\lCap_\CCC(Z) = \llim_{\eps\to 0} \frac{ \log \Lambda_\CCC(Z,\eps)}{\log(1/\eta(\eps))}, \qquad
\uCap_\CCC(Z) = \ulim_{\eps\to 0} \frac{ \log \Lambda_\CCC(Z,\eps)}{\log(1/\eta(\eps))},
\]
where as before the limits are taken over those values of $\eps>0$ for which $\tilde\PPP(Z,\eps)\neq\emptyset$.

The quantities $\dim_\CCC(Z)$, $\lCap_\CCC(Z)$, and $\uCap_\CCC(Z)$ satisfy the following properties.  (For proofs, see~\cite{yP98}.)
\begin{enumerate}
\item $\dim_\CCC \emptyset \leq 0$.
\item $\dim_\CCC Z_1 \leq \dim_\CCC Z_2$ if $Z_1\subset Z_2 \subset X$.
\item $\dim_\CCC \left(\bigcup_{i\geq 0} Z_i \right) = \sup_{i\geq 0} \dim_\CCC Z_i$, where $Z_i \subset X$, $i=0,1,2,\dots$
\item $\dim_\CCC Z \leq \lCap_\CCC Z \leq \uCap_\CCC Z$ for any $Z\subset X$.
\item $\lCap_\CCC Z_1 \leq \lCap_\CCC Z_2$ and $\uCap_\CCC Z_1 \leq \uCap_\CCC Z_2$ for any $Z_1 \subset Z_2 \subset X$.
\item For any sets $Z_i \subset X$, $i=1,2,\dots$,
\[
 \lCap_\CCC \left( \bigcup_{i\geq 0} Z_i \right) \geq \sup_{i\geq 0} \lCap_\CCC Z_i, \qquad
\uCap_\CCC \left( \bigcup_{i\geq 0} Z_i \right) \geq \sup_{i\geq 0} \uCap_\CCC Z_i.
\]
\end{enumerate}

Many dimensional quantities defined for dynamical systems, such as topological pressure and $u$-dimension, are actually given not by a single C-structure, but by a family of C-structures.  That is, one considers for each $\delta>0$ a C-structure $\CCC_\delta$ and then defines the Carath\'eodory dimension by
\[
\dim_\CCC(Z) := \lim_{\delta\to 0} \dim_{\CCC_\delta}(Z),
\]
provided the limit exists.  One may also consider C-structures indexed by a directed set (such as the family of finite open covers of $X$) and take a limit with respect to the ordering on that set.

\begin{proposition}\label{prop:udim}
If $u\in C(X)$ satisfies \PP, then the $u$-dimension defined in Section \ref{sec:pressure} is a Carath\'eodory dimension characteristic.
\end{proposition}
\begin{proof}
Fix $\delta>0$ and let the index set $\SSS$ be the union of the index sets from the covers $E_N$ in \PP.  Let $\FFF$ be the corresponding Bowen balls $B(x,n,\delta)$, with $\eta(x,n) = e^{-S_n u(x)}$, $\xi(x,n) = 1$, and $\psi(x,n) = 1/n$.  Including $\emptyset$ in $\FFF$, this satisfies A1, and A3 follows since each $E_N$ is a cover of $X$ on which $\psi(x,n) = 1/n \leq 1/N$.  Finally, A2 follows from the requirement that $\lim_{N\to\infty} \inf_{(x,n)\in E_N} S_n u(x) = \infty$.
\end{proof}

\bibliographystyle{alpha}
\bibliography{abbrev,books,thermodynamic,multifractal}

\end{document}

%% file: defs.tex
\newcommand{\RR}{\mathbb{R}}

\newcommand{\NN}{\mathbb{N}}

\newcommand{\AAA}{\mathcal{A}}
\newcommand{\BBB}{\mathcal{B}}
\newcommand{\CCC}{\mathcal{C}}

\newcommand{\MMM}{\mathcal{M}}
\newcommand{\HHH}{\mathcal{H}}

\newcommand{\SSS}{\mathcal{S}}
\newcommand{\GGG}{\mathcal{G}}

\newcommand{\ph}{\varphi}
\newcommand{\eps}{\varepsilon}
\newcommand{\di}{\partial}

\newcommand{\llim}{\varliminf}
\newcommand{\ulim}{\varlimsup}
\newcommand{\lr}{\underline{r}}
\newcommand{\ur}{\overline{r}}

\newcommand{\ld}{\underline{d}}
\newcommand{\ud}{\overline{d}}
\newcommand{\htop}{h_\mathrm{top}\,}

\DeclareMathOperator{\diam}{diam}